\documentclass[11pt,reqno,a4paper,oneside]{amsart}

\usepackage[margin=1.75cm]{geometry}
\usepackage{latexsym}
\usepackage{amsmath}
\usepackage{amsthm}
\usepackage{amssymb}
\usepackage{graphicx}
\usepackage{caption}
\usepackage{subcaption}
\usepackage{amsfonts}
\usepackage{xcolor}
\usepackage{fancyhdr}
\usepackage{comment}
\usepackage[section]{placeins}
\usepackage{enumerate}
\usepackage{enumitem}

\usepackage{hyperref}
\usepackage{mathrsfs}                           

\newcommand{\bK}{\mathbf{K}}                    

\newcommand{\eps}{\varepsilon}                
\newcommand{\bx}{\mathbf{x}}
\newcommand{\be}{\mathbf{e}}
\newcommand{\by}{\mathbf{y}}
\newcommand{\bF}{\mathbf{F}}

\newcommand{\bzero}{\mathbf{0}}

\newcommand{\bs}{\mathbf{s}}
\newcommand{\mS}{\mathcal{S}}

\newcommand{\mF}{\mathcal{F}}
\newcommand{\mE}{\mathcal{E}}

\newcommand{\ba}{\mathbf{a}}
\newcommand{\bn}{\mathbf{n}}

\newcommand{\bP}{\mathbf{P}}
\newcommand{\bW}{\mathbf{W}}
\newcommand{\bff}{\mathbf{f}}
\newcommand{\bg}{\mathbf{g}}
\newcommand{\mI}{\mathcal{I}}
\newcommand{\mG}{\mathcal{G}}
\newcommand{\mK}{\mathcal{K}}
\newcommand{\mN}{\mathcal{N}}

\newcommand{\BE}{\begin{equation}}

\newcommand{\bxi}{\boldsymbol\xi}
\newcommand{\bchi}{\boldsymbol\chi}
\newcommand{\bkappa}{\boldsymbol\kappa}
\newcommand{\bes}{\begin{subequations}}
\newcommand{\EE}{\end{equation}}
\newcommand{\ees}{\end{subequations}}
\newcommand{\mO}{\mathcal{O}}
\newcommand{\mM}{\mathcal{M}}

\newcommand{\mQ}{\mathcal{Q}}
\newcommand{\bE}{\mathbf{E}}

\newcommand{\mH}{\mathcal{H}}

\newcommand{\noI}{\noindent}
\newcommand{\BEU}{\begin{equation*}}
\newcommand{\EEU}{\end{equation*}}
\newcommand{\red}[1]{\textcolor{black}{#1}}

\newcommand\mo{
  \mathchoice
    {{\scriptstyle\mathcal{O}}}
    {{\scriptstyle\mathcal{O}}}
    {{\scriptscriptstyle\mathcal{O}}}
    {\scalebox{.7}{$\scriptscriptstyle\mathcal{O}$}}
  }

\numberwithin{equation}{section}

\title[Oscillatory translational instabilities in general 2-D domains]{Oscillatory translational instabilities of localized spot patterns in the Schnakenberg reaction-diffusion system on general 2-D domains}
\author{J.~C.~Tzou$^1$\,,\enspace S. Xie$^{2,*}$
   }

\thanks{$^1$School of Mathematical and Physical Sciences, Macquarie University, Sydney, NSW, Australia; {\tt tzou.justin@gmail.com}}
			
\thanks{$^2$School of Mathematics, Hunan University, Changsha City, Hunan, China; {\tt xieshuangquan2013@gmail.com}}

\thanks{$^*$Corresponding author}

\begin{document}
\begin{abstract}

For a bounded 2-D planar domain $\Omega$, we investigate the impact of domain geometry on oscillatory translational instabilities of $N$-spot equilibrium solutions for a singularly perturbed Schnakenberg reaction-diffusion system with $\mO(\eps^2) \ll \mO(1)$ activator-inhibitor diffusivity ratio. An $N$-spot equilibrium is characterized by an activator concentration that is exponentially small everywhere in $\Omega$ except in $N$ well-separated localized regions of $\mO(\eps)$ extent. We use the method of matched asymptotic analysis to analyze Hopf bifurcation thresholds above which the equilibrium becomes unstable to translational perturbations, which result in $\mO(\eps^2)$-frequency oscillations in the locations of the spots. We find that stability to these perturbations is governed by a nonlinear matrix-eigenvalue problem, the eigenvector of which is a $2N$-vector that characterizes the possible modes (directions) of oscillation. The $2N\times 2N$ matrix contains terms associated with a certain Green's function on $\Omega$, which encodes geometric effects. For the special case of a perturbed disk with radius in polar coordinates $r = 1 + \sigma f(\theta)$ with \red{$0< \varepsilon \ll \sigma \ll 1$}, $\theta \in [0,2\pi)$, and $f(\theta)$ $2\pi$-periodic, we show that only the mode-$2$ coefficients of the Fourier series of $f$ impact the bifurcation threshold at leading order in $\sigma$. We further show that when $f(\theta) = \cos2\theta$, the dominant mode of oscillation is in the direction parallel to the longer axis of the perturbed disk. Numerical investigations on the full Schnakenberg PDE are performed for various domains $\Omega$ and $N$-spot equilibria to confirm asymptotic results and also to demonstrate how domain geometry impacts thresholds and dominant oscillation modes. 

\end{abstract}
\maketitle
\section{Introduction}

The stability and dynamics of spatially localized spike and spot patterns in activator-inhibitor reaction-diffusion systems has been the subject of many studies. These patterns deviate significantly from the uniform state and arise in parameter regimes well beyond the Turing stability threshold, and cannot be well-described by amplitude equations obtained from weakly nonlinear theory \cite{matkowsky1970nonlinear, borckmans1995turing, de1999spatial}. Motivated largely by the 1998 review article \cite{ni1998diffusion}, numerous studies have focused on the so-called semistrong interaction regime \cite{doelman2003semistrong} in which the activator-inhibitor diffusivity ratio $\eps^2 \ll 1$ is asymptotically small. Early works developed matched asymptotic and geometric singular perturbation techniques to characterize the existence, dynamics, and stability of localized patterns in this regime \cite{doelman1998stability, doelman2001large, doelman2000slowly, doelman2003semistrong, iron2001stability, ward2003hopf} of the Gray-Scott and Gierer-Meinhardt systems.

Slow drift dynamics of quasi-equilibrium spot patterns have been determined both asymptotically and numerically for one-, two-, and three-dimensions (see, e.g., \cite{doelman2000slowly, ward2002dynamics, sherratt2013pattern, bastiaansen2019dynamics, carter2018traveling, kolokolnikov2009spot, tzou2017stability, tzou2019spot, bastiaansen2020pulse, wong2021spot}). It is \red{shown} in such works how the drift dynamics are impacted by combinations of advection, interaction with domain boundaries and other spots, domain heterogeneities, and/or curvature, 

Instabilities of spot patterns can occur both in the spot profiles (amplitude instabilities) as well as in the spot locations (translational instabilities). Monotonic (in time) amplitude instabilities come in the form of competition or overcrowding instabilities leading to spot annihilation \cite{doelman2001slowly, sun2005slow, chen2011stability}, or self-replication instabilities \cite{doelman2001slowly, doelman1998stability,muratov2001spike, kolokolnikov2005existence, kolokolnikov2009spot} leading to additional spots. For \red{the} former, it has recently been shown for the 1-D Gierer-Meinhardt and Schnakenberg systems that such competition instabilities are subcritical \cite{kolokolnikov2021competition}. Monotonic translational instabilities of equilibrium spot configurations lead to a rearrangement of spot locations, and may trigger amplitude instabilities in the process. Recently, slow monotonic translations instabilities were analyzed in \cite{kolokolnikov2022ring} for spots in ring equilibrium configurations of the two-dimensional Schnakenberg model. Analysis of these instabilities are often more intricate than that for amplitude instabilities due to the asymptotically small eigenvalues associated with translation instabilities \cite{WGM1, iron2001stability, ward2002existence}.

Results for oscillatory amplitude instabilities of one-dimensional spot patterns have been established for various reaction-diffusion systems; see, e.g., \cite{doelman2002stability, doelman2001large, ward2003hopf, kolokolnikov2005existence, chen2009oscillatory, doelman2015explicit, tzou2, tzou2019spot}. This instability is typically associated with a pair of complex eigenvalues that remain $\mO(1)$ as $\eps \to 0$. In \cite{veerman2015breathing}, a weakly nonlinear theory is developed to characterize oscillations beyond the linear stability regime and determine whether the Hopf bifurcation is subcritical or supercritical. In two dimensions, \cite{wei2001pattern, WGM1,WGS2, WGS1, WS1} established the existence of stability to oscillatory amplitude instabilities, while \cite{tzou2018anomalous} determined an anomalous scaling for the Hopf stability threshold. Recently, Hopf bifurcations for amplitude instabilities were determined for the three-dimensional Gierer-Meinhardt model \cite{gomez2021asymptotic}.

Hopf bifurcations leading to temporal oscillations in spot locations, particularly in dimensions greater than one, has not be analyzed as in-depth as oscillatory amplitude instabilities. In one dimension, \cite{chen2009oscillatory, wei2006slow} obtained Hopf stability thresholds for oscillatory drift instabilities of spots in the Gray-Scott model, and established the $\mO(\eps^2)$ scaling of the associated eigenvalue. Oscillatory instabilities in spot widths, referred to as breathing pulses, were analyzed for three-component systems in \cite{doelman2009pulse, gurevich2013moving, gurevich2006breathing}. \red{In \cite{xie2021complex}, oscillatory motion of multiple spots were investigated for an extended three-component Schnakenberg model, where multiple modes were excited slightly beyond the Hopf bifurcation.}  In two-dimensions, \cite{xie2017moving} determines  oscillatory translational instabilities of a one-spot equilibrium of the Schnakenberg  reaction-diffusion system \cite{schnakenberg1979simple} on the unit disk. The symmetry of the disk, however, meant that the results shed very little light on effects of domain geometry. Furthermore, as the analysis was specific to a \red{one-spot} pattern on the unit disk, it did not account for effects that arise from \red{spot interactions}.


In this paper, we perform the analysis on a general bounded 2-D planar domain $\Omega$ and analyze how its geometry impacts the properties of the instability. In particular, we demonstrate that asymmetries of the domain shape lead to preferred directions of oscillation at the onset of instability, which then saturates into orbits that are far from circular. This is in contrast to the behavior observed in \cite{xie2017moving} for the unit disk. There, it was shown that the Hopf bifurcation of an equilibrium spot at the center of a unit disk was not associated with a preferred direction of motion, and subsequent spot trajectories about the center were nearly circular. We also further generalize the analysis of \cite{xie2017moving} to the case of multi-spot equilibria, and deduce how spot orientation affects the mode of oscillations.

For simplicity, we will consider the now well-studied (dimensionless) Schnakenberg reaction-diffusion model  for activator $v(\bx,t)$ and inhibitor $u(\bx,t)$ concentrations

\bes \label{schnak}
\BE \label{schnakv}
	v_t = \eps^2\Delta v - v + uv^2 \,,  \qquad t > 0 \,, \quad \bx = \left(\begin{array}{c} x_1 \\ x_2 \end{array}\right) \in \Omega \,,
\EE
\BE \label{schnaku}
	\tau u_t = \Delta u + A - \eps^{-2} uv^2 \,,  \qquad t > 0 \,, \quad \bx \in \Omega \,,
\EE
\BE \label{schnakbc}
	\partial_n v = 0 \enspace \mbox{and} \enspace \partial_n \red{u}= 0 \,, \quad \bx \in \partial\Omega  \,,
\EE
\ees

where $\partial_n$ denotes the normal derivative on $\partial \Omega$. We note that any non-unity diffusion coefficient $D$ on the inhibitor component can be scaled out independently of domain size. That is, with a diffusivity $D$ for the inhibitor, we recover  \eqref{schnaku} by setting $v \to \sqrt{D}v$, $u \to u/\sqrt{D}$, $\tau \to D\tau$, $A \to \sqrt{D} A$.

In \eqref{schnak}, $\eps^2 \ll 1$ is the small diffusivity of the activator, $A$ is a (constant) feed rate, and the Hopf bifurcation parameter $\tau$ is a measure of the rate at which the inhibitor responds to perturbations in the concentrations of the activator and inhibitor. As $\tau$ is increased, an increasingly sluggish response of the inhibitor leads to oscillatory instabilities via Hopf bifurcations \cite{kerner2013autosolitons}. 

The outline of the paper is as follows. In \S \ref{Nspoteqbm}, we asymptotically construct an $N$-spot equilibrium solution of the Schnakenberg PDE\eqref{schnak}. In contrast to the $\mO(1)$ construction presented in \cite{kolokolnikov2009spot}, we require correction terms up to $\mO(\eps^2)$ in order to facilitate the subsequent stability analysis. In \S \ref{stab} we analyze the stability of the $N$-spot equilibrium to oscillatory translation instabilities. We derive a $2N\times 2N$ complex matrix-eigenvalue problem of the form $P\mathbf{a} = \lambda\mathbf{a}$ that characterizes the Hopf bifurcation of a translational perturbation of a one-spot equilibrium. The Hopf bifurcation threshold for $\tau$ is obtained by requiring that the $2N\times 2N$ matrix $P$ have a pure imaginary eigenvalue $\lambda$, which yields the frequency of oscillations. The eigenvector $\ba$ will yield initial directions along which spot oscillations occur. Effects of domain geometry are encoded in the entries of $P$, which involve the quadratic terms of the local behavior of Helmholtz Green's function $G_\mu(\bx;\bx_0)$ satisfying (see \cite{paquin2020asymptotics} for a derivation of \eqref{Gkloc})

\bes \label{Gk}
\BE \label{Gkeq}
	\Delta G_\mu - \mu G_\mu = -\delta(\bx-\bx_j) \,, \quad \bx,\bx_j \in \Omega \,, \qquad \partial_n G_\mu = 0 \,, \quad \bx \in \partial \Omega \,,
\EE
\begin{multline} \label{Gkloc}
	G_\mu \sim -\frac{1}{2\pi}\log|\bx-\bx_j| + R_\mu(\bx_j;\bx_j)  + \nabla_\bx R_\mu (\bx;\bx_j)\left.\right|_{\bx=\bx_j} \cdot (\bx-\bx_j) \\  - \frac{\mu}{8\pi}|\bx-\bx_j|^2 \log|\bx-\bx_j|  + \frac{1}{2}(\bx-\bx_j)^T H_{\mu j}(\bx-\bx_j) \,, \enspace \mbox{as} \enspace \bx\to\bx_j \,,
\end{multline}
and 
\begin{multline} \label{Gklocij}
	G_\mu \sim G_\mu(\bx_i;\bx_j) + \nabla_{\bx}G_\mu(\bx;\bx_j)\left.\right|_{\bx=\bx_j} \cdot (\bx-\bx_j) + \frac{1}{2}(\bx-\bx_i)H_{\mu ij}(\bx-\bx_i) \,, \enspace \mbox{as} \enspace \bx\to\bx_i \,,
\end{multline}

\ees

where $H_{\mu j}$ and $H_{\mu ij}$ are the respective $2\times 2$ Hessian matrices. Their entries depend on the geometry of the domain (and also on $\mu$) and thus cannot be obtained from a local analysis. With the exception of special geometries such as disks and rectangles, they must be computed numerically.

In \S \ref{onespot}, we reduce our $2N\times 2N$ eigenvalue problem result to the case of $N=1$. We compare the resulting eigenvalue problem to the one derived for the special case of one spot inside the unit disk in \cite{xie2017moving}. We use this comparison to highlight extra terms in the analysis that arise due to asymmetries of domain geometry. In \S \ref{perturbeddisk}, we analyze how perturbations of the unit disk break the symmetry and give rise to two distinct thresholds corresponding to two distinct modes of oscillation. We show that the lower of these thresholds, and therefore, the preferred mode of oscillation, is determined by the mode-2 coefficient of the Fourier series of the perturbation. In \S \ref{numerics}, we perform detailed numerical investigations to confirm our theoretical results for both the $1$- and $N$-spot cases by solving the full Schnakenberg PDE \eqref{schnak}. In these computations, we consider domains that are high in symmetry (e.g., half disk, unit disk, rectangles) and domains that have little symmetry (e.g., rectangular domains containing circular holes).



\section{Equilibrium and stability analysis}
%
%
%

In this section, we investigate the impact of domain geometry on the preferred initial direction of oscillation of the oscillatory translational instability of $N$-spot equilibrium solutions to \eqref{schnak}. We begin with a brief construction of the equilibrium solutions before performing the stability analysis.

\subsection{$N$-spot equilibrium} \label{Nspoteqbm}

For completeness, we begin with a very brief outline of the construction of an $N$-spot equilibrium; more details can be found in, e.g., \cite{kolokolnikov2009spot}. In the inner region near the $j$-th spot centered at $\bx = \bx_j$, we have the inner variables

\bes\label{innervar} 
\BE \label{innervar1} 
	\bx = \bx_j + \eps\by_j \,, \quad  \quad \by_j = \left(\begin{array}{c} y_{1j} \\ y_{2j} \end{array}\right) =  \rho_j \be_j \,, \quad \be_j \equiv \left(\begin{array}{c} \cos \theta_j \\ \sin\theta_j \end{array}\right) \,,
\EE
\BE \label{innervar2} 
	u_e\sim U_{0j}(\rho_j) + \eps^2 U_{2j}(\by_j)\,, \quad v_e \sim V_{0j}(\rho_j) + \eps^2 V_{2j}(\by_j) \,,
\EE
\ees

where the $\mO(\eps)$ terms in \eqref{innervar2} are absent under under the assumption that each spot is stationary in time. Note also that the leading order spot profile is radially symmetric; the asymmetry due to the geometry of the problem is captured at $\mO(\eps^2)$. Substituting \eqref{innervar} into \eqref{schnak} and collecting leading order terms, we obtain the following core problem for the radially symmetric functions $U_{0j}$ and $V_{0j}$,

\bes \label{core}
\BE \label{coreeq}
	\Delta_{\rho_j} V_{0j} - V_{0j} + U_{0j}V_{0j}^2 = 0 \,, \quad \Delta_{\rho_j} U_{0j} - U_{0j}V_{0j}^2 = 0 \,, \qquad \rho_j > 0
\EE
\BE \label{corebc}
	V_{0j}^\prime(0) = U_{0j}^\prime(0) = 0 \,, \quad V_{0j} \to 0 \enspace \mbox{and} \enspace U_{0j} \sim S_j\log\rho_j + \chi(S_j) \,, \enspace \mbox{as} \enspace \rho_j \to \infty \,.
\EE
\ees

In \eqref{coreeq}, $\Delta_{\rho_j} \equiv \partial_{\rho_j\rho_j} + \rho_j^{-1}\partial_{\rho_j}$ denotes the radially symmetric Laplacian in the polar coordinates $(\rho_j, \theta_j)$. Numerical solutions of \eqref{core} are depicted in Fig. 2 of \cite{kolokolnikov2009spot} for various $S_j$, including the nonlinear function $\chi(S_j)$. We assume that $S_j \lessapprox 4.3$ so that each spot is stable to the local (mode-2) self-replication instability. The divergence theorem on \eqref{core} yields

\BE \label{Sdiv}
	2\pi S_j = \int_{\mathbb{R}^2}\! U_{0j} V_{0j}^2 \, d\by \,.
\EE 

%
%

In the outer region, the $\eps^{-2}uv^2$ term in \eqref{schnaku} behaves in the distributional sense as a sum of delta functions located at each $\bx_j$ with weight $2\pi S_j$ as given in \eqref{Sdiv} As such, the leading order equilibrium solution $u_e(\bx)$ satisfies 

\BE \label{ueeq}
	\Delta u_e + A = 2\pi \sum_{j=1}^N S_j \delta(\bx-\bx_j) \,, \qquad \bx \in \Omega \,; \quad \partial_n u_e = 0 \enspace \mbox{on} \enspace \bx \in \Omega \,,
\EE

with solution given by

\BE \label{uout}
	u_e(\bx) \sim -2\pi \sum_{j = 1}^N S_j G(\bx;\bx_j) + \bar{u} \,,
\EE

where by the zero-integral condition on $G$ below, the constant $\bar{u}$ is the average of $u$ over $\Omega$. Integrating \eqref{ueeq} over $\Omega$ and applying the divergence theorem yields the solvability condition

\BE \label{Sj}
	2\pi\sum_{j = 1}^N S_j = A|\Omega| \,,
\EE

where $|\Omega|$ denotes the area of the domain $\Omega$. In \eqref{uout}, $G(\bx;\bxi)$ is the modified Neumann Green's function satisfying

\bes \label{GNall}
\BE \label{GNeq}
	\Delta G(\bx;\bx_j) = -\delta(\bx-\bx_j) + \frac{1}{|\Omega|} \,, \quad \bx,\bx_j \in \Omega \,, \qquad \partial_n G = 0 \,, \quad \bx \in \partial\Omega \,; \qquad \int_\Omega \! G (\bx;\bx_j) \,d\bx = 0 \,,
\EE

\BE \label{GNloc}
	G \sim -\frac{1}{2\pi}\log|\bx-\bx_j| + R_{jj} + \nabla R_{jj} \cdot (\bx-\bx_j) +  \frac{1}{2}(\bx-\bx_j)^T H_{jj}(\bx-\bx_j) \enspace \mbox{as} \enspace \bx \to \bx_j \,,
\EE

where $R_{jj}$ is the regular part of $G(\bx;\bx_j)$ evaluated on the diagonal, $\nabla R_{jj} = \nabla_\bx R(\bx;\bx_j) \left.\right|_{\bx=\bx_j}$ is its gradient, and $H_{jj}$ its Hessian matrix.  For $\bx_i \neq \bx_j$, we have that

\BE \label{Gij}
	G(\bx;\bx_i) \sim G_{ij} + \nabla G_{ji} \cdot (\bx-\bx_i) + \frac{1}{2} (\bx-\bx_i)^T H_{ji}(\bx-\bx_i) \,, \enspace \mbox{as} \enspace \bx \to \bx_j \,,
\EE

where $G_{ji} = G(\bx_j;\bx_i)$, while $\nabla G_{ji}$ and $H_{ji}$ are the gradient and Hessian terms of $G(\bx;\bx_i)$ at $\bx_j$, respectively, in the Taylor expansion. Principal Result 3.4 in \cite{kolokolnikov2009spot} gives the equation of motion $d\bx_j/dt$ of the $j$-th spot in terms of the gradient terms $\nabla R_{jj}$ and $\nabla G_{ji}$ with $i \neq j$. The condition $d\bx_j/dt = \bzero$ for all $j = 1,\ldots,N$ yields the $2N$ equations for the equilibrium strengths $S_j$ and locations $\bx_j$
\ees

\BE \label{eqbm1}
	S_j\nabla R_{jj} + \sum_{i\neq j}^N S_i \nabla G_{ji} = \bzero \,, \quad j = 1,\ldots, N \,.
\EE

The equations \eqref{Sj} and \eqref{eqbm1} constitute $2N+1$ equations for $S_j$, $\bx_j$, and $\bar{u}$. To determine the remaining $N$ equations that fix the parameters of an $N$-spot equilibrium, we match the far-field behavior of the $j$-th inner region \eqref{corebc} to the leading order terms of the local behavior of the outer solution \eqref{uout} near $\bx_j$, yielding the $N$ nonlinear equations

\bes \label{matchjj}
\BE \label{matchj} 
\left(2\pi \mathcal{G}\nu + \mathcal{I}_N\right) \mathbf{s} + \nu\bchi  = \nu\bar{u} \mathbf{e} \,;
\EE
with
\BE \label{matchj2} 
\mathbf{G} \equiv  \left(\begin{array}{cccc} R_{11} & G_{12} & \cdots & G_{1N} \\  G_{21} & \ddots & \ddots & \vdots \\  \vdots & \ddots & \ddots & G_{N-1, N} \\  G_{N1} & \hdots & G_{N,N-1}  & G_{NN} \end{array} \right) \,, \quad \mathbf{s} = \left(\begin{array}{c} S_1 \\ \vdots \\ S_N \end{array}\right) \,, \quad \mathbf{e} = \left(\begin{array}{c} 1 \\ \vdots \\ 1 \end{array}\right) \,, \quad \bchi = \left(\begin{array}{c} \chi(S_1)\\ \vdots \\ \chi(S_N)  \end{array}\right) \,,
\EE
\ees

where $\mI_m$ is the $m\times m$ identity matrix. Equation \eqref{matchjj} along with \eqref{Sj} and \eqref{eqbm1} determine the spot strengths $S_j$, the spot locations $\bx_j$, along with $\bar{u}$ that define a leading order equilibrium $N$-spot configuration.

The $\mO(\eps^2)$ correction terms $U_{2j}(\by_j)$ and $V_{2j}(\by_j)$ satisfy the system for $\by_j \in \mathbb{R}^2$,

\bes
\BE \label{coreeq2}
	\Delta_{\by_j} V_{2j} - V_{2j} + 2U_{0j}V_{0j}V_{2j} + V_{0j}^2U_{2j} = 0 \,, \quad \Delta_{\by_j} U_{2j} - 2U_{0j}V_{0j}V_{2j} - V_{0j}^2U_{2j} = 0 \,, \qquad \by_j \in \mathbb{R}^2 \,.
\EE

The far-field conditions for $U_{2j}$ and $V_{2j}$ come from the quadratic terms in the local behavior of \red{$u_e$} near $\bx_j$. That is, from  $u_e$ in \eqref{uout} and local behaviors of $G(\bx;\bx_j)$ in \eqref{GNall}, as $\bx \to \bx_j$, 

\begin{multline*}
u_e \sim -2\pi \left[ -\frac{S_j}{2\pi} \log|\bx-\bx_j| + S_j R_{jj} + S_j\nabla R_jj \cdot (\bx-\bx_j) + \frac{1}{2}S_j(\bx-\bx_j)^TH_{jj}(\bx-\bx_j)  \right. \\ \left.  + \sum_{i\neq j} G_{ji} S_i + \sum_{i\neq j} S_i \nabla G_{ji} \cdot (\bx-\bx_j) + \frac{1}{2} \sum_{i\neq j} S_i (\bx-\bx_j)^T H_{ji} (\bx-\bx_j) \right] + \bar{u} \,.
\end{multline*}

All terms not involving $(\bx-\bx_j)$ are matched at $\mO(1)$, yielding \eqref{matchjj}. Terms linear in $(\bx-\bx_j)$ sum to zero due to the equilibrium condition \eqref{eqbm1}, which is why the inner expansion \eqref{innervar2} has no $\mO(\eps)$ term, while the quadratic terms are matched by the far-field of $U_{2j}$. With the inner polar coordinates $\rho_j$ and $\be_j$ as defined in \eqref{innervar1}, we have
 the far-field conditions 

\BE \label{innervar2far}
	U_{2j} \sim -\pi\rho_j^2 \be_j^T \mH_j \be_j \,, \qquad V_{2j} \to 0 \,, \enspace \mbox{as} \enspace \rho_j \to \infty \,,
\EE

where the $2\times 2$ matrix $\mH_j$ is defined as

\BE \label{mHj}
	\mH_j \equiv \sum_{i=1}^N S_i H_{ji} \,.
\EE

\ees

This completes the construction of the $N$-spot equilibrium solution to \eqref{schnak}, with the inner solution given to $\mO(\eps^2)$ by \eqref{innervar} along with \eqref{core}, \eqref{coreeq2}, and \eqref{innervar2far}. In the outer region, the leading order equilibrium for $u$ is given in \eqref{uout}, while $v_e \sim 0$.

\subsection{Stability analysis} \label{stab} 

With $|\psi|,|\phi| \ll \mO(1)$, we perturb the equilibrium solution

\BE \label{pertmain}
	u \sim u_e + e^{\lambda t} \psi + c.c. \,, \qquad v \sim v_e + e^{\lambda t} \phi + c.c. \,,
\EE

which yields the eigenvalue problem

\BE \label{eigprob}
	\lambda\phi = \eps^2 \Delta\phi - \phi + 2u_ev_e\phi + v_e^2\psi \,, \quad \tau \lambda \psi = \Delta \psi -\frac{1}{\eps^2}\left[ 2u_ev_e\phi + v_e^2\psi  \right] \,.
\EE

In \eqref{pertmain},  $c.c.$ denotes the complex conjugate of the term immediately preceding it.

In the inner region, we let $\psi \sim \Psi_j(\by_j)$ and $\phi \sim \Phi_j(\by_j)$, where 

\BE \label{pert}
	\Psi_j \sim \Psi_{0j} + \eps\Psi_{1j} + \eps^2\Psi_{2j} \,, \qquad \Phi_j \sim \Phi_{0j} + \eps\Phi_{1j} + \eps^2\Phi_{2j} \,.
\EE

Since drift velocities of spots in quasi-equilibrium patterns are $\mO(\eps^2)$ (see e.g., \cite{kolokolnikov2009spot}), we assume that $\lambda \sim \mO(\eps^2)$ and that $\tau\lambda \sim \mO(1)$ when $\tau$ is at or near the Hopf bifurcation threshold. Substituting \eqref{pert} into \eqref{schnak} and collecting leading order terms, we have for $\Psi_{0j}$ and $\Phi_{0j}$

\bes \label{pertsys}
\BE \label{O1pert}
	\Delta_{\by_j} \left(\begin{array}{c} \Phi_{0j} \\ \Psi_{0j}   \end{array}\right) + \mM_j\left(\begin{array}{c} \Phi_{0j} \\ \Psi_{0j}   \end{array}\right) = \bzero \,.
\EE

where 

\BE \label{M}
	\mM_j \equiv \left(\begin{array}{cc} -1 + 2U_{0j}V_{0j} & V_{0j}^2 \\ -2U_{0j}V_{0j} & - V_{0j}^2   \end{array}\right) \,.
\EE

We observe that any linear combination of $\partial_{y_{1j}} V_{0j}$ and $\partial_{y_{2j}} V_{0j}$ satisfies the first equation in \eqref{O1pert}, while any linear combination of $\partial_{y_{1j}} U_{0j}$ and $\partial_{y_{2j}} U_{0j}$ satisfies the second. That is, the perturbation is the translation mode given by

\BE \label{O1pertsol}
	\Phi_{0j} = \ba_j^T\nabla_{\by_j} V_{0j} = \partial_{\rho_j} V_{0j} (\ba_j^T\be_j) \,, \qquad \Psi_{0j} = \ba_j^T\nabla_{\by_j} U_{0j} = \partial_{\rho_j} U_{0j} (\ba_j^T\be_j) \,; \qquad \ba_j \equiv \left(\begin{array}{c} a_{1j} \\ a_{2j} \end{array}\right) \,.
\EE

In \eqref{O1pertsol}, the possibly complex vector $\ba_j$ determines the nature of the oscillations of the $j$-th spot at the onset of instability. If $\ba_j$ is real, the $j$-th spot oscillates along a line passing through the equilibrium location $\bx_j$ in the direction $\ba_j$; if $\ba_j$ is complex, the motion of the $j$-th spot is rotational about the point $\bx_j$. From \eqref{corebc} and \eqref{O1pertsol}, the far-field behaviors of $\Phi_{0j}$ and $\Psi_{0j}$ are 

\BE \label{O1pertfar}
\Phi_{0j} \to 0 \,, \qquad \Psi_{0j} \sim (\ba_j^T\be_j)\frac{S_j}{\rho_j} \,, \enspace \mbox{as} \enspace \rho_j \to \infty \,.
\EE
\ees

We contrast \eqref{pertsys} with stability analysis of amplitude instabilities; i.e., the self-replication peanut instability (e.g., \cite{kolokolnikov2009spot}, competition instability (e.g., \cite{chen2011stability}), and amplitude oscillation instabilities (e.g., \cite{tzou2018anomalous}). These amplitude instabilities all occur on an $\mO(1)$ time-scale so that $\lambda \sim \mO(1)$, giving rise to a $\lambda \Phi_{0j}$ in the right-hand side of the equation for $\Phi_{0j}$ in \eqref{O1pert}. Furthermore, the competition and amplitude oscillation instabilities are radially symmetric to leading order with $\Psi_{0j} \sim \log\rho_j$ in the far-field, which leads to a strong coupling through the inhibitor component in the outer region. On the other hand, the self-replication eigenfunction is a mode-2 instability with $\Psi_{0j} \sim \rho_j^{-2}$ in the far-field. The fast decay leads to a very weak coupling between the different spots. It is therefore a strictly local instability, to leading order. The translation instability is mode-1 with $\Psi_{0j} \sim \rho_j^{-1}$. This decay leads to a weak coupling through an $\mO(\eps)$ outer solution for the inhibitor eigenfunction. In contrast to the local self-replication instability, we show below that the nature of this weak coupling between the spots must be determined in order to characterize the translation instability.

The $1/\rho_j$ behavior of $\Psi_{0j}$ in the far-field gives rise to a singular behavior in the outer region near $\bx_j$ that must be matched by the leading order term of $\psi$, which is $\mO(\eps)$. The regular part of the behavior of $\psi$ at $\bx_j$ must then be matched by a constant term in the far-field of $\Psi_{1j}$. At $\mO(\eps)$ in the inner region, we have for $\Phi_{1j}$ and $\Psi_{1j}$,

\bes \label{Oepspert}
\BE \label{Oepsperteq}
	\Delta_{\by_j} \left(\begin{array}{c} \Phi_{1j} \\ \Psi_{1j}   \end{array}\right) + \mM_j\left(\begin{array}{c} \Phi_{1j} \\ \Psi_{1j}   \end{array}\right) = \bzero \,.
\EE

Since $\Psi_{1j}$ must have a constant term in the far-field, we must have

\BE \label{Oepspertfar}
\Phi_{1j} \to 0 \,, \qquad \Psi_{1j} \sim \kappa_j(\nu) \left[\log\rho_j + B_j(S_j)   \right] \,, \enspace \mbox{as} \enspace \rho_j \to \infty \,.
\EE
\ees

In \eqref{Oepspertfar}, $\kappa_j(\nu)$ is a scaling constant to be found by matching to the regular part of $\psi_1$ at the spot locations. The solution to \eqref{Oepspert} is (see e.g., \cite{tzou2017stability})

\BE \label{Phi1Psi1sol}
	\Phi_{1j} = \kappa_j\partial_{S_j} V_{0j} \,, \qquad \Psi_{1j} = \kappa_j\partial_{S_j} U_{0j} \,,
\EE

where $\partial_{S_j} U_{0j} \sim \log\rho_j + \chi^\prime(S_j)$, and $\chi(S_j)$ is the constant that must be computed from the core problem \eqref{core}. The far-field behavior of $\Psi_{1j}$ is thus

\BE \label{Psi1jfar}
	\Psi_{1j} \sim \kappa_j \left[\log|\by_j| + \chi^\prime_j \right] \,, \enspace \mbox{as} \enspace |\by_j| \to \infty \,,
\EE

where $\chi^\prime_j \equiv \chi^\prime(S_j)$. Observe that both the dipole term from $\Psi_{0j}$ as well as the logarithmic term of $\Psi_{1j}$ must be contained in the singularity structure for $\psi$ in the outer region.

Thus, in the outer region, we let $\psi \sim \eps \psi_1$, where $\psi_1$ satisfies

\bes \label{psi1}
\BE \label{psi1eq}
	\Delta \psi_1 = \tau\lambda\psi_1 \,, \qquad \bx \in \Omega \,, \qquad \partial_n \psi_1 = 0 \,, \quad \bx \in \partial\Omega
\EE

with the singular behavior 

\BE \label{psi1sing}
\psi_1 \sim S_j \frac{\ba_j^T(\bx-\bx_j)}{|\bx-\bx_j|^2} + \kappa_j\log|\bx-\bx_j| + \kappa_j\left[\frac{1}{\nu} + \chi^\prime_j \right] \,, \enspace \mbox{as} \enspace \bx \to \bx_j \,, \quad j = 1, \ldots, N \,.
\EE
\ees

In terms of the Helmholtz Green's function of \eqref{Gk} and its gradient with respect to the second variable, we determine $\psi_1$ to be

\BE \label{psi1sol}
 \psi_1 = 2\pi\sum_{i = 1}^N \left[ S_i \ba_i^T \nabla_{\bx_i} G_{\lambda\tau}(\bx;\bx_i) - \kappa_i G_{\lambda\tau}(\bx;\bx_i) \right] \,.
\EE

Notice that $\nabla_{\bx_i} G_{\lambda\tau}(\bx;\bx_i)$ produces the dipole term of \eqref{psi1sing} near $\bx_i$ while it also still satisfies the no-flux boundary condition of \eqref{psi1eq} since the gradient is being taken with respect to the second variable. Its local behaviors near $\bx_i$ and $\bx_j \neq \bx_i$ are

\bes \label{gradGk}
\begin{multline} \label{gradGKii}
\nabla_{\bx_i} G_{\lambda\tau}(\bx;\bx_i) \sim \frac{1}{2\pi}\frac{\bx-\bx_i}{|\bx-\bx_i|^2} + \bF_{{\lambda\tau} i} + \mF_{{\lambda\tau} i}(\bx-\bx_i) \\ + \frac{{\lambda\tau}}{4\pi}(\bx-\bx_i)\log|\bx-\bx_i| + \left[\frac{{\lambda\tau} }{8\pi} \mI_2 - H_{{\lambda\tau}_i} \right](\bx-\bx_i) \,, \enspace \mbox{as} \enspace \bx \to \bx_i \,;
\end{multline}
\BE \label{gradGKji}
\nabla_{\bx_i} G_{\lambda\tau}(\bx;\bx_i) \sim \bE_{{\lambda\tau} ji} + \mE_{{\lambda\tau} ji}(\bx-\bx_j) \ \,, \enspace \mbox{as} \enspace \bx \to \bx_j \neq \bx_i \,;
\EE

where we have defined the quantities

\begin{equation} \label{gradGKdef}
\begin{gathered}
	\bF_{\mu_i} = \left(\begin{array}{c} F_{\mu_i}^{(1)}  \\ F_{\mu_i}^{(2)}   \end{array}\right) \equiv  \nabla_\bx R_\mu(\bx;\bx_i)\left.\right|_{\bx=\bx_i} \,, \qquad \mF_{\mu_i} \equiv \left(\nabla_{\bx_i} F_{\mu_i}^{(1)}   \enspace\enspace \nabla_{\bx_i} F_{\mu_i}^{(2)}   \right) \,,  \\
	\bE_{\mu_{ji}} \equiv \nabla_{\bx_i} G_\mu(\bx_j;\bx_i) \,, \qquad \mE_{\mu ji} \equiv \left(\nabla_{\bx_i}\partial_{x_1} G_\mu(\bx;\bx_i) \left.\right|_{\bx = \bx_j}   \enspace\enspace \nabla_{\bx_i}\partial_{x_2} G_\mu(\bx;\bx_i) \left.\right|_{\bx = \bx_j}  \right) \,.
\end{gathered}
\end{equation}
\ees

The scaling constant $\kappa_j$ of \eqref{Oepspertfar} is then found by matching the far-field of $\Psi_{1j}$ in \eqref{Psi1jfar} to the constant terms of the local behavior of $\psi_1$ near $\bx_j$. Using \eqref{psi1sol}, \eqref{gradGk}, and \eqref{Gk}, we match the constant terms in \eqref{psi1sing} and those contained in \eqref{psi1sol} near $\bx_j$ to  obtain the matching condition

\BE  \label{kappaj}
	\kappa_j \left[\frac{1}{\nu} + \chi_j^\prime + 2\pi R_{\lambda\tau}(\bx_j;\bx_j)\right] + 2\pi \sum_{i\neq j}\kappa_i G_{\lambda\tau_{ji}} = 2\pi S_j\ba_j^T \bF_{\lambda\tau_j} + 2\pi \sum_{i \neq j} S_{i} \ba_i^T \bE_{\lambda\tau_{ji}} \,, \quad j = 1, \ldots, N \,.
\EE

In matrix-vector form, \eqref{kappaj} becomes

\bes \label{kappajvec}
\BE \label{kappajeq}
	\bkappa = \mK_{\lambda\tau} \ba
\EE
where we have defined the $N\times 1$ vector $\bkappa$, $2N\times 1$ vector $\ba$, the $N\times N$ matrices $\mG_{\lambda\tau}$ and $\bchi^\prime$, the $N\times 2N$ matrices $(\nabla_2 \mG_{\lambda\tau})$ and $\mK_{\lambda\tau}$, and the diagonal $2N \times 2N$ matrix $\mS$
\BE \label{kappadef}
\begin{gathered}
	\bkappa \equiv \left(\begin{array}{c}  \kappa_1 \\ \vdots \\ \kappa_N\end{array}\right) \,, \quad \ba \equiv \left(\begin{array}{c}  \ba_1 \\ \vdots \\ \ba_N \end{array}\right) \,, \quad   \mG_{\lambda\tau} \equiv \left(\begin{array}{cccc} R_{\lambda\tau}(\bx_1;\bx_1) & G_{\lambda\tau}(\bx_1;\bx_2) & \cdots & G_{\lambda\tau}(\bx_1;\bx_N) \\ G_{\lambda\tau}(\bx_2;\bx_1) & R_{\lambda\tau}(\bx_2;\bx_2) & \cdots & \vdots \\ \vdots & \cdots & \ddots & \vdots \\ G_{\lambda\tau}(\bx_N;\bx_1) & \cdots & \cdots & R_{\lambda\tau}(\bx_N;\bx_N) \end{array}\right)   \\ \bchi^\prime \equiv \left(\begin{array}{ccc} \chi_1^\prime & & \\ & \ddots & \\ & & \chi_N^\prime \end{array}\right) \,, \quad (\nabla_2 \mG_{\lambda\tau}) \equiv \left(\begin{array}{cccc} \bF_{\lambda\tau_1}^T & \bE_{\lambda\tau_{12}}^T & \cdots & \bE_{\lambda\tau_{1N}}^T \\ \bE_{\lambda\tau_{21}}^T & \bF_{\lambda\tau_2}^T & \cdots & \vdots \\ \vdots & \cdots & \ddots & \vdots \\ \bE_{\lambda\tau_{N1}}^T & \cdots & \cdots & \bF_{\lambda\tau_N}^T \end{array}\right) \,, \quad \mS \equiv \left(\begin{array}{ccccc} S_1 & & & & \\  & S_1 & & & \\   & & \ddots & & \\  & & & S_N & \\   & & & & S_N \end{array}\right) \,, \\ \mK_{\lambda\tau} \equiv  2\pi \left[ \frac{1}{\nu} \mI_N + \bchi^\prime + 2\pi \mG_{\lambda\tau}\right]^{-1} (\nabla_2 \mG_{\lambda\tau}) \mS \,.
\end{gathered}
\EE
\ees

In \eqref{kappadef}, the vectors $\bE_{\mu ij}$ and $\bF_{\mu j}$  are gradients of the Green's function and its regular part, respectively, defined in \eqref{gradGKdef}. Also, we note that as $\nu \to 0$, the matrix to be inverted in the computation of $\mK_{\lambda\tau}$ in \eqref{kappadef} is invertible due to its being diagonally dominant. 

The linear system for $\kappa_j$ in \eqref{kappajvec} along with \eqref{psi1sol} and  \eqref{Phi1Psi1sol} determine the leading order outer solution for $\psi$, up to the oscillation directions $\ba_j$, and the $\mO(\eps)$ inner solutions for $\phi$ and $\psi$. The Hopf stability threshold for $\tau$, the frequency of oscillations at onset $\lambda$, and the direction of oscillations $\ba_j$ will be determined via a solvability condition at $\mO(\eps^2)$ for $\Phi_{2j}$ and $\Psi_{2j}$. To proceed, we must first obtain the far-field condition for $\Psi_{2j}$ from the linear and $|\bx-\bx_j|\log|\bx-\bx_j|$ terms in the local behavior of $\psi_1$ near $\bx_j$. Recalling that $\psi_1$ is an $\mO(\eps)$ term, while $\bx-\bx_j = \eps\rho_j\be_j$, we use \eqref{psi1sol} with \eqref{Gk} and \eqref{gradGk} to compute that as $\rho_j \to \infty$,

\begin{multline} \label{Psi2jfar1}
\Psi_{2j} \sim \frac{1}{2} S_j\lambda\tau (\ba_j^T\be_j)\rho_j\log\rho_j + 2\pi\rho_j \left\lbrace S_j\ba_j^T\left[\mF_{\lambda\tau_j} - H_{\lambda\tau j} + \frac{\lambda\tau}{8\pi}\left(1-\frac{2}{\nu}   \right)\mI_2\right] - \right.  \\ \left.  \kappa_j \bF_{\lambda\tau_j}^T - \sum_{i \neq j}\left[ \kappa_i \nabla_{\bx}G_{\lambda\tau}(\bx;\bx_i)\mid_{\bx=\bx_j}^T + S_i\ba_i^T \mE_{\lambda\tau_{ji}}\right]   \right\rbrace \be_j \,, \quad j = 1,\ldots,N \,.
\end{multline}

Next, we substitute \eqref{pert} into \eqref{eigprob} while recalling the expansion for $u_e$ and $v_e$ in \eqref{innervar2}, and with $\lambda = \eps^2\lambda_0$ and $\tau = \eps^{-2}\tau_0$, we collect $\mO(\eps^2)$ terms to obtain

\BE \label{Oeps2eq}
\begin{gathered}
	\lambda_0\Phi_{0j} = \Delta_{\by_j} \Phi_{2j} - \Phi_{2j} + 2U_{0j}V_{0j} \Phi_{2j} + V_{0j}^2 \Psi_{2j} + 2\left(U_{0j}V_{2j} + U_{2j} V_{0j}\right) \Phi_{0j} + 2U_{0j}V_{2j} \Psi_{0j} \,, \\
	\lambda_0\tau_0 \Psi_{2j} = \Delta_{\by_j} \Psi_{2j} - 2U_{0j}V_{0j} \Phi_{2j} - V_{0j}^2 \Psi_{2j} - 2\left(U_{0j}V_{2j} + U_{2j} V_{0j}\right) \Phi_{0j} - 2U_{0j}V_{2j} \Psi_{0j} \,.
\end{gathered}
\EE

To write \eqref{Oeps2eq} and \eqref{Psi2jfar1} more compactly, we define the quantities

\BE \label{Oeps2def}
\begin{gathered}
	\bW_j \equiv \left(\begin{array}{c} \Phi_{2j} \\ \Psi_{2j}   \end{array}\right) \,, \qquad \bff_{1j} \equiv \left(\begin{array}{c} \Phi_{0j} \\ \Psi_{0j}   \end{array}\right) = \ba_j^T\be_j \left(\begin{array}{c} \partial_{\rho_j} V_{0j} \\ \partial_{\rho_j} U_{0j}   \end{array}\right) \,, \qquad \omega \equiv \lambda_0\tau_0 \\ \bff_{2j} \equiv \left(\begin{array}{c} \lambda_0 \Phi_{0j} \\ \omega \Psi_{0j}   \end{array}\right) = \ba_j^T\be_j \left(\begin{array}{c} \lambda_0 \partial_{\rho_j} V_{0j} \\ \omega \partial_{\rho_j} U_{0j}   \end{array}\right) \,, \qquad \mN_j \equiv \left(\begin{array}{cc} -2(U_{0j}V_{2j} + U_{2j}V_{0j}) & -2V_{0j}V_{2j} \\  2(U_{0j}V_{2j} + U_{2j}V_{0j}) & 2V_{0j}V_{2j}   \end{array}\right) \,.
\end{gathered}
\EE

where we have used \eqref{O1pertsol} for $\Phi_{0j}$ and $\Psi_{0j}$, to obtain

\bes \label{Oeps2inn}
\BE \label{Oeps2inneq}
 \Delta_{\by_j} \bW_j + \mM_j \bW_j = \mN_j \bff_{1j} + \bff_{2j} \,, \qquad \by_j \in \mathbb{R}^2
\EE

with the far-field condition

\BE \label{Psi2jfar2}
	\mathbf{W}_j \sim \left(\begin{array}{c} 0 \\ \left[\frac{1}{2} S_j\omega\log\rho_j \ba_j^T + 2\pi\ba_{Q_j}^T\right]\rho_j\be_j \end{array}\right) \,, \enspace \mbox{as} \enspace \rho_j \to \infty \,, \quad j = 1\, \ldots, N \,, 
\EE

where we have defined the $N\times 2$ matrix $(\nabla_1\mG_\omega)_j$ and $2N\times 2$ matrices $M_{\omega j}$ and $Q_{\omega j}$, and $2\times 1$ vector $\ba_{Qj}$,

\BE \label{LjMj}
\begin{gathered}
 (\nabla_1\mG_\omega)_j \equiv \left(\begin{array}{c} \nabla_{\bx}G_{\omega}(\bx;\bx_1)\mid_{\bx=\bx_j}^T  \\ \nabla_{\bx}G_{\omega}(\bx;\bx_2)\mid_{\bx=\bx_j}^T \\ \vdots \\ \bF_{\omega j}^T \\ \vdots \\ \nabla_{\bx}G_{\omega}(\bx;\bx_N)\mid_{\bx=\bx_j}^T   \end{array} \right)	\,, \qquad  M_{\omega_j} \equiv \left(\begin{array}{c} S_1 \mE_{\omega_{j1}}  \\ S_2\mE_{\omega_{j2}} \\ \vdots \\ S_j\left[\mF_{\omega_j} - H_{\omega_j} + \frac{\omega}{8\pi}\left(1-\frac{2}{\nu}\right)\mI_2 \right] \\ \vdots \\ S_N\mE_{\omega_{jN}}    \end{array} \right) \,, \\ Q_{\omega_j} \equiv M_{\omega_j} - \mK_\omega^T (\nabla_1\mG_\omega)_j \,, \qquad \ba_{Qj} \equiv \left(\begin{array}{c} a_{Q_{1j}} \\ a_{Q_{2j}} \end{array}\right) =  Q_{\omega j}^T \ba \,.
\end{gathered}
\EE

\ees

In \eqref{Oeps2inneq} and \eqref{Psi2jfar2} for $\bW_j$, the coupling of the $j$-th inner region is contained only in the $\ba_{Qj}$ term defined in \eqref{LjMj}. All other terms are local to the $j$-th inner region.

From \eqref{pertsys}, the linear operator in \eqref{Oeps2inneq} admits a nontrivial nullspace of dimension at least two. The nonhomogeneous terms of \eqref{Oeps2inneq} and \eqref{Psi2jfar2} must therefore satisfy an orthogonality condition involving the solution to the homogeneous adjoint problem. Before applying this condition, we observe that $\bW_j$ can be decomposed into two components proportional to $\cos\theta_j$ and $\sin\theta_j$, respectively, as

\bes \label{bW}
\BE \label{bWsol}
\bW_j = \bW_{cj}\cos\theta_j + \bW_{sj}\sin\theta_j \,,
\EE

where $\bW_{cj}$ and $\bW_{sj}$ are $2\times 1$ vectors satisfying

\BE  \label{bWeq}
\begin{gathered}
	\left(\partial_{\rho_j\rho_j} + \frac{1}{\rho_j}\partial_{\rho_j} - \frac{1}{\rho_j^2}   \right) \bW_{cj} + \mM_j\bW_{cj} = a_{1j}\left[\mN_j \left(\begin{array}{c} \partial_{\rho_j} V_{0j} \\ \partial_{\rho_j} U_{0j} \end{array}\right) +  \left(\begin{array}{c} \lambda_0 \partial_{\rho_j} V_{0j}  \\ \omega \partial_{\rho_j} U_{0j} \end{array}\right) \right] \,, \\ \left(\partial_{\rho_j\rho_j} + \frac{1}{\rho_j}\partial_{\rho_j} - \frac{1}{\rho_j^2}   \right) \bW_{sj} + \mM_j\bW_{sj} = a_{2j}\left[\mN_j \left(\begin{array}{c} \partial_{\rho_j} V_{0j} \\ \partial_{\rho_j} U_{0j} \end{array}\right) +  \left(\begin{array}{c} \lambda_0 \partial_{\rho_j} V_{0j} \\ \omega \partial_{\rho_j} U_{0j} \end{array}\right) \right] \,,
\end{gathered}
\EE

with the far-field conditions

\BE  \label{bWfar}
\begin{gathered}
	\bW_{cj} \sim \left(\begin{array}{c} 0 \\ \frac{1}{2}S_j\omega a_{1j}\rho_j\log\rho_j + 2\pi\rho_j a_{Q_{1j}} \end{array}\right) \,, \enspace \mbox{as} \enspace \rho_j \to \infty\,.\\ 	\bW_{sj} \sim \left(\begin{array}{c} 0 \\ \frac{1}{2}S_j\omega a_{2j}\rho_j\log\rho_j + 2\pi\rho_j a_{Q_{2j}} \end{array}\right) \,, \enspace \mbox{as} \enspace \rho_j \to \infty \,.
\end{gathered}
\EE

\ees

The nonhomogeneous terms of \eqref{Oeps2inn} must be orthogonal to the nullspace of the homogeneous adjoint operator, given by

\BE \label{adjeq}
\Delta_{\by_j} \bP_j + \mM_j^T \bP_j = \bzero \,, \qquad \by_j \in \mathbb{R}^2 \,.
\EE

We seek two linearly independent mode-1 solutions To \eqref{adjeq} of the form $\bP_j = \bP_{cj}\cos\theta_j$ and $\bP_j = \bP_{sj}\sin\theta_j$, where $\bP_{cj}$ and $\bP_{sj}$ are given by

\BE \label{adjdecomp}
\bP_{cj} \equiv \tilde{\bP}_j(\rho_j)\cos\theta_j \,, \qquad \bP_{sj} \equiv \tilde{\bP}_j(\rho_j)\sin\theta_j \,, \qquad \tilde{\bP}_j(\rho_j) \equiv \left(\begin{array}{c} \tilde{P}_{1j} \\ \tilde{P}_{2j} \end{array}\right)
\EE

and the radially symmetric $\tilde{\bP}_j$ satisfies

\bes \label{Ptilde}
\BE \label{Ptildeeq}
	\left(\partial_{\rho_j\rho_j} + \frac{1}{\rho_j}\partial_{\rho_j} - \frac{1}{\rho_j^2}   \right) \tilde{\bP}_j + \mM_j^T\tilde{\bP}_j = \bzero\,, \qquad 0 < \rho_j < \infty
\EE

with boundary and far-field conditions

\BE \label{Ptildefar}
	\tilde{\bP}_j(\bzero) = \bzero \,, \qquad \tilde{\bP}_j \sim \left(\begin{array}{c} 0 \\ 1/\rho_j \end{array}\right)  \,, \enspace \mbox{as} \enspace \rho_j \to \infty \,.
\EE
\ees

Note that the normalization condition in the far-field condition of \eqref{Ptildefar} uniquely specifies $\tilde{\bP}_j$, while the condition at the origin ensures continuity of $\bP_{cj}$ and $\bP_{sj}$. To apply the orthogonality condition, we multiply \eqref{Oeps2inn} on the left by $\bP_{cj,sj}^T$ and integrate over a disk $B_{R}$ of radius $R \gg 1 $ centered at the origin to obtain

\BE \label{orthog}
\iint_{B_R} \! \bP_{cj,sj}^T \left[\Delta_{\by_j} \bW_j + \mM_j \bW_j \right] \, d\by_j \ = \iint_{B_R} \! \bP_{cj,sj}^T \mN_j \bff_{1j} \, d\by_j + \iint_{B_R} \! \bP_{c,s}^T \bff_{2j} \, d\by_j \,.
\EE 

We now compute each of the integrals in \eqref{orthog} in the limit $R \gg 1$. For the term on the left-hand side, we use Green's identity along with \eqref{adjeq} to obtain

\BE \label{orthogleft}
\iint_{B_R} \! \bP_{cj,sj}^T \left[\Delta_{\by_j} \bW_j + \mM_j \bW_j \right] \, d\by_j = \int_{\theta_j = 0}^{2\pi}\! \left[\bP_{cj,sj}^T \partial_{\rho_j}\bW_j - \bW_j^T\partial_{\rho_j} \bP_{cj,sj}   \right] R \, d\theta_j \,.
\EE \,.

For the right-hand side of \eqref{orthogleft}, we use the far-field conditions of \eqref{bWfar} along with \eqref{adjdecomp} and \eqref{Ptildefar} to obtain for the cosine term

\bes \label{orthogleftcossin}
\BE \label{orthogleftcos}
\int_{\theta_j = 0}^{2\pi}\! \left[\bP_{cj}^T \partial_{\rho_j}\bW_j - \bW_j^T\partial_{\rho_j} \bP_{cj}   \right] R \, d\theta_j \sim  \pi\left[2c_{1j}\log R + 2c_{2j} + c_{1j} \right] \,, \qquad R \gg 1 \,,
\EE

while for sine term, we have

\BE\label{orthogleftsin}
\int_{\theta_j = 0}^{2\pi}\! \left[\bP_{sj}^T \partial_{\rho_j}\bW_j - \bW_j^T\partial_{\rho_j} \bP_{sj}   \right] R \, d\theta_j \sim  \pi\left[2s_{1j}\log R + 2s_{2j} + s_{1j} \right] \,, \qquad R \gg 1 \,,
\EE

where we have defined

\BE \label{orthogleftdef}
c_{1j} \equiv \frac{1}{2}S_j\omega a_{j1} \,, \qquad c_{2j} \equiv 2\pi a_{Qj1}  \,,  \qquad s_{1j} = \frac{1}{2}S_j\omega a_{j2} \,, \qquad s_{2j} \equiv 2\pi a_{Qj2} \,.
\EE

\ees

For the second term on the right-hand side of \eqref{orthog}, we use \eqref{Oeps2def} for $\bff_{2j}$ and perform an integration by parts to obtain

\BE \label{orthogPf2}
\begin{gathered}
	\iint_{B_R} \! \bP_{cj}^T \bff_{2j} \, d\by_j \sim \pi a_{1j}\omega S_j\log R + \pi a_{1j}\lambda_0\left[-\tau_0 k_{2j} + k_{1j}\right] \,,\\
	\iint_{B_R} \! \bP_{sj}^T \bff_{2j} \, d\by_j \sim \pi a_{2j}\omega S_j\log R + \pi a_{2j}\lambda_0\left[-\tau_0 k_{2j} + k_{1j}\right] \,,
\end{gathered}
\EE

where $k_{1j}$ and $k_{2j}$ are defined by the integrals (see \cite{xie2017moving})

\BE \label{k1jk2j}
	k_{1j} \equiv \int_0^\infty \! V_{0j}^\prime \tilde{P}_{1j} \rho_j \, d\rho_j \,, \qquad k_{2j} \equiv \int_0^\infty \! \left[U_{0j} - \chi_j\right] (\tilde{P}_{2j} \rho_j)^\prime \, d\rho_j \,.
\EE

\begin{figure}[!ht]
     \centering
     \begin{subfigure}[b]{0.48\textwidth}
         \centering
         \includegraphics[width=\textwidth]{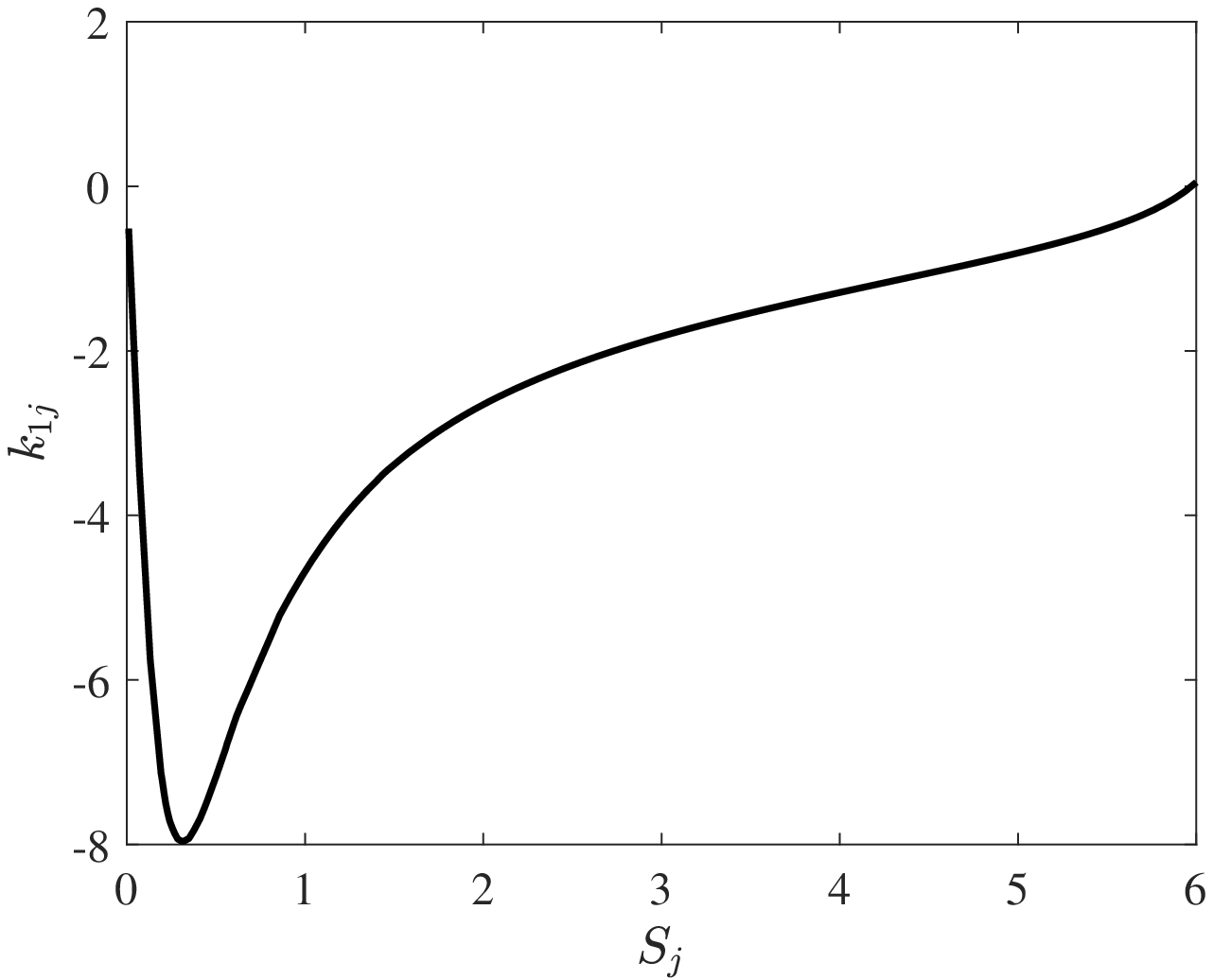}
         \caption{$k_{1j}$}
         \label{kappa_1}
     \end{subfigure}
     \begin{subfigure}[b]{0.48\textwidth}
         \centering
         \includegraphics[width=\textwidth]{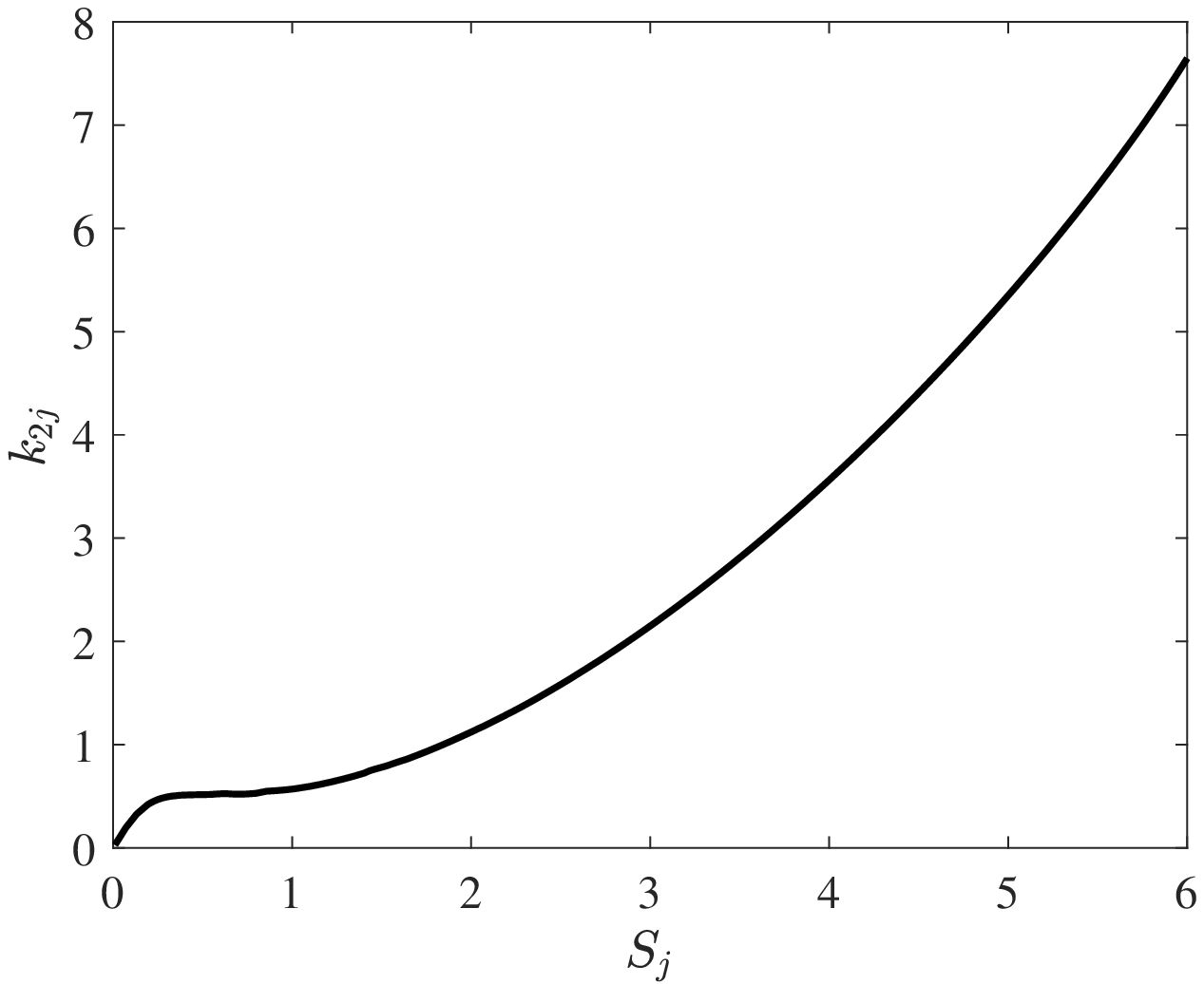}
         \caption{$k_{2j}$}
         \label{kappa_2 }
     \end{subfigure}
        \caption{Plots of (a) $k_{1j}$ and (b) $k_{2j}$ as defined in \eqref{k1jk2j}. Note that a spot with strength $S_j \gtrsim 4.31$ is unstable to a self-replication instability.}
        \label{kappafig}
\end{figure}

Note that $k_{1j}$ and $k_{2j}$ are both functions of $S_j$ through their dependence on $V_{0_j}$ and $\tilde{P}_{1j}$, as well as $U_{0_j}$, $\tilde{P}_{2j}$, and $\chi_j$, respectively. For completeness, we reproduce plots of $k_{1j}$ and $k_{2j}$ versus $S_j$ in Fig \ref{kappafig}. 

\renewcommand{\thesubfigure}{\alph{subfigure}}

To compute the first term on the right-hand side of \eqref{orthog}, we seek to rewrite $\mN_j \bff_1$ in terms of the linear operator $\Delta_{\by_j}  + \mM_j$ so that, upon multiplication by $\bP_{cj,sj}^T$ and integrating over $B_{R}$, we may apply the divergence theorem. To proceed, we use \eqref{Oeps2def} to compute

\BE \label{Nf1_1}
	\mN_j \bff_1 = \mN_j (a_{1j} \cos\theta_j\partial_{\rho_j} + a_{2j}\sin\theta_j\partial_{\rho_j})  \left(\begin{array}{c} V_{0j} \\  U_{0j} \end{array}\right) = \mN_j \ba_j \cdot \nabla_{\by_j}\left(\begin{array}{c} V_{0j} \\  U_{0j} \end{array}\right) \,.
\EE

Expanding \eqref{Nf1_1}, we have

\BE \label{Nf1_2}
	\mN_j \bff_1 = \left(\begin{array}{c}   U_{2j}\ba_j\cdot\nabla_{\by_j} V_{0j}^2 + 2V_{2j}\ba_j\cdot\nabla_{\by_j}(U_{0j}V_{0j})     \\   -U_{2j}\ba_j\cdot\nabla_{\by_j} V_{0j}^2 - 2V_{2j}\ba_j\cdot\nabla_{\by_j}(U_{0j}V_{0j})    \end{array}\right) \,.
\EE

Next, we add and subtract terms to each component of \eqref{Nf1_2} to obtain full derivative derivatives of $U_{2j}V_{0j}^2$ and $2U_{0j}V_{0j}V_{2j}$ and find

\BE \label{Nf1_3}
	\mN_j \bff_1 = \left(\begin{array}{c}   \ba_j\cdot\nabla_{\by_j}(U_{2j} V_{0j}^2) + \ba_j\cdot\nabla_{\by_j}(2V_{2j}U_{0j}V_{0j}) - V_{0j}^2 \ba_j\cdot\nabla_{\by_j} U_{2j} - 2U_{0j}V_{0j}\ba_j\cdot\nabla_{\by_j} V_{2j}     \\   -\ba_j\cdot\nabla_{\by_j}(U_{2j} V_{0j}^2) - \ba_j\cdot\nabla_{\by_j}(2V_{2j}U_{0j}V_{0j}) + V_{0j}^2 \ba_j\cdot\nabla_{\by_j} U_{2j} + 2U_{0j}V_{0j}\ba_j\cdot\nabla_{\by_j} V_{2j}      \end{array}\right) \,.
\EE

Passing the operator $\ba_j\cdot\nabla_{\by_j}$ through the system of $U_{2j}$ and $V_{2j}$ in \eqref{coreeq2}, we observe that 

\BEU
	\ba_j\cdot\nabla_{\by_j}(U_{2j} V_{0j}^2) + \ba_j\cdot\nabla_{\by_j}(2V_{2j}U_{0j}V_{0j}) = -\Delta_{\by_j} \ba_j\cdot\nabla_{\by_j} V_{2j} + \ba_j\cdot\nabla_{\by_j} V_{2j} \,,
\EEU

and

\BEU
	\ba_j\cdot\nabla_{\by_j}(U_{2j} V_{0j}^2) + \ba_j\cdot\nabla_{\by_j}(2V_{2j}U_{0j}V_{0j}) = \Delta_{\by_j} \ba_j\cdot\nabla_{\by_j} U_{2j} \,.
\EEU

We thus obtain from \eqref{Nf1_3} that

\BE \label{Nf1_4}
	\mN_j \bff_1 = \left[\Delta_{\by_j} + \mM_j\right]\ba_j\cdot\nabla_{\by_j} \bn_j \,; \qquad \bn_j \equiv \left(\begin{array}{c}  V_{2j}  \\   U_{2j}   \end{array}\right) \,,
\EE

where $\mM_j$ is the matrix of the linearized reaction terms defined in \eqref{M}. Multiplying \eqref{Nf1_4} by $\bP_{cj,sj}^T$ satisfying \eqref{adjeq} and applying the divergence theorem, we have

\BE \label{PNf1_1}
 \iint_{B_R} \! \bP_{cj,sj}^T \mN_j \bff_{1j} \, d\by_j = \int_0^{2\pi}\left(\bP_{cj,sj}^T \partial_{\rho_j}(\ba_j\cdot\nabla_{\by_j} \bn_j) - (\ba_j\cdot\nabla_{\by_j} \bn_j)^T\partial_{\rho_j} \bP_{cj,sj}\right)R\,d\theta_j \,,
\EE

where each term of the integrand on the right-hand side are evaluated on the circle $\rho_j = R$. For $\bn_j$ defined in \eqref{Nf1_4}, we use the far-field condition for $V_{2j}$ and $U_{2j}$ in \eqref{innervar2far} to compute 

\BE \label{njfar}
\begin{gathered}
\ba_j\cdot\nabla_{\by_j} \bn_j \sim \left(\begin{array}{c} 0\\ -2\pi\rho_j  \ba_j^T\mH_j \be_j \end{array}\right) \,; \qquad \rho_j \gg 1 \,, \\
\partial_{\rho_j} \ba_j\cdot\nabla_{\by_j} \bn_j \sim  \left(\begin{array}{c} 0\\ -2\pi\ba_j^T\mH_j \be_j \end{array}\right) \,; \qquad \rho_j \gg 1 \,.
\end{gathered}
\EE

Substituting \eqref{njfar} and \eqref{adjdecomp}-\eqref{Ptildefar} into \eqref{PNf1_1}, we obtain

\BE \label{PNf1_2}
\begin{gathered}
\iint_{B_R} \! \bP_{cj}^T \mN_j \bff_{1j} \, d\by_j \sim -4\pi^2\left[a_{j1}\mH_{j}^{(11)} + a_{j2}\mH_{j}^{(21)}\right] \,,   \\
\iint_{B_R} \! \bP_{sj}^T \mN_j \bff_{1j} \, d\by_j \sim -4\pi^2\left[a_{j1}\mH_{j}^{(12)} + a_{j2}\mH_{j}^{(22)}\right] \,, 
\end{gathered}
\EE

where $\mH_j^{(mn)}$ denotes the $(m,n)$-th entry of the matrix $\mH_j$ defined in \eqref{mHj}.

Now we substitute \eqref{orthogleftcossin}, \eqref{orthogPf2} and \eqref{PNf1_2} for their respective terms in the orthogonality condition \eqref{orthog} to obtain for $R \gg 1$ and $j = 1,\ldots,N$,

\BE \label{orthog2}
\begin{gathered}
	\pi\left[2c_{1j}\log R + 2c_{2j} + c_{1j} \right] = \pi a_{j1}\omega S_j\log R + \pi a_{j1}\lambda_0\left[-\tau_0 k_{2j} + k_{1j}\right] - 4\pi^2\left[a_{j1}\mH_{j}^{(11)} + a_{j2}\mH_{j}^{(21)}\right]  \,, \\
	\pi\left[2s_{1j}\log R + 2s_{2j} + s_{1j} \right] = \pi a_{j2}\omega S_j\log R + \pi a_{j2}\lambda_0\left[-\tau_0 k_{2j} + k_{1j}\right] -4\pi^2\left[a_{j1}\mH_{j}^{(12)} + a_{j2}\mH_{j}^{(22)}\right] \,. 
\end{gathered}
\EE

Observe that, with $c_{1j}$ and $s_{1j}$ defined in terms of $a_{j1}$ and $a_{j2}$ in \eqref{orthogleftdef}, the $\log R$ terms in \eqref{orthog2} cancel, yielding

\BE \label{mateig}
\begin{gathered}
	\left[\frac{1}{2}S_j\omega + \omega k_{2j} + 4\pi^2\mH_j^{(11)}\right] a_{j1} + 4\pi^2\mH_j^{(21)}a_{j2} + 4\pi a_{Qj1} = k_{1j}\lambda_0a_{j1} \,, \\
	4\pi^2\mH_j^{(12)} a_{j1} + \left[\frac{1}{2}S_j\omega + \omega k_{2j} + 4\pi^2\mH_j^{(22)}\right]a_{j2} + 4\pi a_{Qj2} = k_{1j}\lambda_0a_{j2} \,.
\end{gathered}
\EE

where $a_{Qj1}$ and $a_{Qj2}$ are the first and second components, respectively, of the vector $\ba_{Qj}$ defined in \eqref{LjMj}, which is the term giving rise to the coupling between the different inner regions. Substituting for $\ba_{Qj}$ in \eqref{mateig}, we arrive at the $2N\times 2N$ matrix-eigenvalue problem for the oscillation modes $\ba$, Hopf bifurcation frequency $\lambda_0$ and threshold $\tau_0 \equiv \omega/\lambda_0$:

\bes \label{finaleig}
\BE \label{finaleigeq}
\bK_1^{-1}\left[\omega\left(\frac{1}{2}\mS + \bK_2\right) + 4\pi\mH + 4\pi \mQ_\omega    \right] \ba = \lambda_0\ba \,.
\EE

where $\mS$ is the $2N\times 2N$ matrix of spot strengths \eqref{kappadef}, while the $2N\times 2N$ matrices $\bK_1$, $\bK_2$, $\mH$, and $\mQ_\omega$ are defined as

\BE \label{finaleigdef}
\begin{gathered}
	\bK_1 \equiv \left(\begin{array}{ccccc} k_{11} & & & & \\  & k_{11} & & & \\   & & \ddots & & \\  & & & k_{1N} & \\   & & & & k_{1N} \end{array}\right) \,, \qquad \bK_2 \equiv \left(\begin{array}{ccccc} k_{21} & & & & \\  & k_{21} & & & \\   & & \ddots & & \\  & & & k_{2N} & \\   & & & & k_{2N} \end{array}\right) \,, \\ \mQ_\omega = \left(\begin{array}{c} Q_{\omega 1}^T \\ Q_{\omega 2}^T \\ \vdots \\ Q_{\omega N}^T  \end{array}\right) \,, \qquad \mH \equiv \left(\begin{array}{cccc} \mH_1  & & & \\  & \mH_2  & & \\   & & \ddots &  \\ &  & & \mH_N    \end{array}\right) \,.
\end{gathered}
\EE
\ees

In \eqref{finaleigdef}, $k_{1j}$ and $k_{2j}$ are the integrals defined in \eqref{k1jk2j}, $Q_{\omega j}$ are the $2N \times 2$ matrices defined in \eqref{LjMj}, and $\mH_{j}$ is the $2\times 2$ matrix defined in \eqref{mHj}. Rewriting in terms of previously defined matrices, we obtain that the final matrix-eigenvalue problem of \eqref{finaleigeq} takes the form

\bes \label{finaleigblock}
\BE 
	\bK_1^{-1}\left[ B(\omega) - M(\omega)  \right] \ba = \lambda_0\ba
\EE

where the $2N\times 2N$ matrices $B(\omega)$ and $M(\omega)$ are given by

\BE \label{finaleigblockdef}
\begin{gathered}
     B(\omega) \equiv \omega\left[\left(1-\frac{1}{\nu}\right)\mS + \bK_2 \right] + 4\pi\mH + 4\pi \left(\nabla^2 \mG_\omega\right)\mS \,, \\ M(\omega) \equiv 8\pi^2\left(\nabla_1 \mG_\omega\right)^T\left[\frac{1}{\nu}\mI_N + \bchi^\prime + 2\pi\mG_\omega\right]^{-1} \left(\nabla_2\mG_\omega\right)\mS \,,
	\end{gathered}
\EE

with the $N\times 2N$ and $2N\times 2N$ matrices $\left(\nabla_1\mG_\omega\right)$ and $\left(\nabla^2\mG_\omega\right)$, respectively, given by

\BE \label{finaleigblockdef2}
\begin{gathered}
 \left(\nabla_1\mG_\omega\right) \equiv \left(\left(\nabla_1\mG_\omega\right)_1 \enspace\enspace \left(\nabla_1\mG_\omega\right)_2 \enspace\enspace  \cdots \enspace\enspace   \left(\nabla_1\mG_\omega\right)_N  \right) \,, \\
(\nabla^2 \mG_{\omega}) \equiv \left(\begin{array}{cccc} [\mF_{\omega_1} - H_{\omega_1}]^T & \mE_{\omega_{12}}^T & \cdots & \mE_{\omega_{1N}}^T \\ \mE_{\omega_{21}}^T & [\mF_{\omega 2} - H_{\omega_2}]^T & \cdots & \vdots \\ \vdots & \cdots & \ddots & \vdots \\ \mE_{\omega_{N1}}^T & \cdots & \cdots & [\mF_{\omega_N} - H_{\omega_N}]^T \end{array}\right) \,.
	\end{gathered}
\EE 
\ees

In \eqref{finaleigblock}, $\mS$ is the diagonal matrix of spot strengths defined in \eqref{kappadef}, $\mH$ is the diagonal block matrix with each block a linear combination of Hessian matrices of the Neumann Green's function weighted by spot strengths (see \eqref{mHj}), $\left(\nabla^2 \mG_\omega\right)$ is a block matrix involving second order derivatives of the Helmholtz Green's function and its regular part (see \eqref{gradGKdef}), $\left(\nabla_1\mG_\omega\right)$ and $\left(\nabla_2\mG_\omega\right)$ are $N\times 2N$ matrices involving first derivatives of the first and second arguments of the Helmholtz Green's function $G_\omega(\bx;\bx_0)$, respectively (see \eqref{finaleigblockdef2} and \eqref{kappadef}), $\mI_N$ is the $N\times N$ identity matrix, $\bchi^\prime$ is the diagonal matrix whose $j$-th diagonal entry is $\chi^\prime(S)$, while $\mG_\omega$ is the Green's interaction matrix  The diagonal matrices $\bK_1$ and $\bK_2$ are defined in \eqref{finaleigdef}; the terms along the diagonal are the nonzero constants defined in \eqref{k1jk2j} that must be computed numerically. The $(i,j)$-th entry or block in each matrix when $i \neq j$ accounts for the interaction between the $i$-th and $j$-th spot.

To find pure imaginary eigenvalues of \eqref{finaleigblock}, we set $\lambda_0 = i\lambda_I$ and solve the two transcendental equations 

\bes \label{det0}
\BE \label{det01}
\textrm{Re}\left\lbrace\det\left(B(i\hat{\omega}) - M(i\hat{\omega}) - i\hat{\lambda}_I\mI_{2N}\right)\right\rbrace = 0 \,,
\EE
\BE \label{det02}
\textrm{Im}\left\lbrace\det\left(B(i\hat{\omega}) - M(i\hat{\omega}) - i\hat{\lambda}_I\mI_{2N}\right)\right\rbrace = 0 \,,
\EE
\ees

for $\omega = \hat{\omega}$ and $\lambda_I = \hat{\lambda}_I$. Each solution corresponds to a different mode of oscillation with frequency $\hat{\lambda}_I/(2\pi)$ and Hopf bifurcation value $\hat{\tau} = \hat{\omega}/\hat{\lambda}_I$. The corresponding nullspace $\hat{\ba}$ of the matrix $B(i\hat{\omega}) - M(i\hat{\omega}) - i\hat{\lambda}_I\mI_{2N}$ then indicates the mode of oscillation of each of the $N$ spots. The Hopf stability threshold $\hat{\tau}^*$ is then the minimum of all $\hat{\tau}$ with the dominant mode of oscillation given by the corresponding eigenvector $\hat{\ba}^*$ and frequency $\lambda_I^*$.

If the direction vector $\hat{\ba}_k$ of the $k$-th spot is real, then at onset, the $k$-th spot oscillates about the point $\bx_k \in \Omega$ along the direction $\hat{\ba}_k$. If $\hat{\ba}_k = \textrm{Re}(\hat{\ba}_k) + i\textrm{Im}(\hat{\ba}_k)$, the trajectory of the $k$-th spot at onset is that of a rotated ellipse centered at $\bx_k$ with angle of rotation and minor and major axes determined by the $2\times 2$ matrix $(\textrm{Re}(\hat{\ba}_k) \enspace -\textrm{Im}(\hat{\ba}_k))$. We give numerical examples of both types of oscillation in \S \ref{numerics}.

\section{Single-spot analysis} \label{onespot}

In this section, we consider the special case of the eigenvalue problem \eqref{finaleigblock} when the pattern consists of only $N=1$ spot. We compare our result to that given in \cite{xie2017moving} for the special case of the unit disk, showing that we recover its result from \eqref{finaleigblock} when $N=1$ and $\Omega$ is the unit disk. In doing so, we highlight the extra terms that arise in the eigenvalue problem due to asymmetries of the domain that were absent in the analysis of the unit disk. We also perform an analysis of a perturbed unit disk and show how the perturbation affects the bifurcation threshold, corresponding frequency, as well as the preferred direction of oscillation.

In the case of a single spot of strength $S$, the only terms that remain in \eqref{finaleigblock} are the (1,1) entries and blocks of the matrices in the formula. Then the small eigenvalue problem for $N=1$ is \eqref{finaleigblock}, where the $2\times 2$ matrices $B(\omega)$ and $M(\omega)$ are given by

\BE \label{oneeigdef}
\begin{gathered}
     B(\omega) \equiv \omega\left[\left(1-\frac{1}{\nu}\right)S + k_{21} \right]\mI_2 + 4\pi H_{11} + 4\pi\left[\mF_{\omega_1} - H_{\omega_1} \right]^T  \,, \\ M(\omega) \equiv 8\pi^2 S \left[\frac{1}{\nu} + \chi^\prime(S) + 2\pi R_\omega(\bx_1;\bx_1)\right]^{-1} \left(\nabla_\bx R_\omega(\bx;\bx_1)\left.\right|_{\bx=\bx_1}\right)\left(\nabla_\bx R_\omega(\bx;\bx_1)\left.\right|_{\bx=\bx_1}\right)^T\,.
	\end{gathered}
\EE

Recall that $\nabla_\bx R(\bx;\bx_1)\mid_{\bx=\bx_1} = \bzero$ since $\bx_1$ is an equilibrium location of the one-spot pattern. In highly symmetric geometries such as rectangles and the unit disk which, the zero of $\nabla_\bx R(\bx;\bx_1)\mid_{\bx=\bx_1}$ coincides exactly with that of $\nabla_\bx R_\omega(\bx;\bx_1)\left.\right|_{\bx=\bx_1}$. That is, the gradient of the regular parts of the Neumann \eqref{GNall} and Helmholtz \eqref{Gk} Green's functions vanish at the same value of $\bx_1$. In these geometries, the matrix $M(\omega)$ vanishes while the Hessian matrices $H_{11}$ and $H_{\omega_1}$ along with $\mF_{\omega_1}$ can be made diagonal. The eigenvalue problem \eqref{finaleigblock} then decouples to form

\BE \label{oneeigsymm1}
\begin{gathered}
	k_{11}^{-1}\left\lbrace \omega\left[\left(1-\frac{1}{\nu}\right)S + k_{21} \right]\mI_2 + 4\pi \left[H_{11}^{(1,1)} + \mF_{\omega_1}^{(1,1)} -  H_{\omega_1}^{(1,1)} \right]S \right\rbrace a_1 = \lambda_0 a_1 \,, \\
	k_{11}^{-1}\left\lbrace \omega\left[\left(1-\frac{1}{\nu}\right)S + k_{21} \right]\mI_2 + 4\pi \left[H_{11}^{(2,2)} + \mF_{\omega_1}^{(2,2)} - H_{\omega_1}^{(2,2)} \right]S \right\rbrace a_2 = \lambda_0 a_2 \,,
	\end{gathered}
\EE

where $H_{11}^{(jj)}$ and $H_{\omega_1}^{(jj)}$ are the $(j,j)$ components of the Hessian matrices of the Neumann and Helmholtz Green's functions, respectively. The dominant mode of oscillation would either along the $(1,0)$ or $(0,1)$ directions depending on which pair of Hessian terms in \eqref{oneeigsymm1} yields the lower Hopf threshold $\hat{\tau}$.

In even more symmetric geometries such as the square and unit disk, the latter of which we analyze in detail below, $H_{11}$ and $H_{\omega_1}$ are multiples of $\mI_2$ so that \eqref{finaleigblock} reduces to the scalar problem given by

\BE \label{oneeigsymm2}
 \omega\left[\left(1-\frac{1}{\nu}\right)S + k_{21} \right] + 4\pi \left[H_{11}^{(1,1)} + \mF_{\omega_1}^{(1,1)} - H_{\omega_1}^{(1,1)} \right]S - k_{11}\lambda_0  = 0 \,,
\EE

In such geometries, the vector $\ba$ indicating the direction of spot oscillation at onset becomes arbitrary, and there is no preferred direction of oscillation. 

Observe from \eqref{kappaj} with $N = 1$ that the coincidence in the zeros of $\nabla_\bx R(\bx;\bx_1)\mid_{\bx=\bx_1}$ and $\nabla_\bx R_\omega(\bx;\bx_1)\left.\right|_{\bx=\bx_1}$ implies that $\kappa_1 = 0$. With $\kappa_1 = 0$, the additional extra $G_{\lambda\tau}$ term in the outer solution for the eigenfunction $\psi_1$ in \eqref{psi1sol} vanishes, while the $\mO(\eps)$ term is also absent in the inner expansion for the eigenfunctions $\Phi_1$ and $\Psi_1$. As such, we may attribute the presence of these two terms to the asymmetry of the domain. Their effects are encoded in the matrix $M(\omega)$ in \eqref{oneeigdef}. 

In the next section, we consider \eqref{oneeigsymm2} for the case of the unit disk and show that it is equivalent to the eigenvalue problem derived in \cite{xie2017moving}.

\subsection{The unit disk} \label{unitdisk}

We first require the Hessian terms of the Neumann and Helmholtz Green's functions, $H_{11}^{(1,1)}$ and $H_{\omega_1}^{(1,1)}$, respectively, along with $\mF_{\omega_1}^{(1,1)}$, the gradient with respect to the source location of the gradient of the regular part of the Helmholtz Green's function.

We begin with computing $\mF_{\omega_1}^{(1,1)}$. In polar coordinates $\bx = (x,y) = \rho(\cos\theta, \sin\theta)$, \cite{chen2011stability} gives the series solution for the Helmholtz Green's function $G_\omega(\bx;\bx_0)$ satisfying \eqref{Gk} with source at $\bx_0 = (x_0, y_0) = \rho_0(\cos\theta_0, \sin\theta_0)$ as

\begin{equation} \label{Gunit}
G_\omega(\rho,\theta;\rho_0,\theta_0) = \frac{1}{2\pi}K_0\left(\sqrt{\omega}|\bx-\bx_0|\right) - \frac{1}{2\pi}A_0I_0(\sqrt{\omega}\rho) - \frac{1}{\pi}\sum_{n  =1}^\infty \cos(n(\theta-\theta_0))A_n I_n(\sqrt{\omega}\rho) \,,
\end{equation}

\noI where $|\bx-\bx_0| = \sqrt{\rho^2 + \rho_0^2 - 2\rho\rho_0\cos(\theta-\theta_0)}$ and $A_n = K_n^\prime(\sqrt{\omega})I_n(\sqrt{\omega}\rho_0)/I_n^\prime(\sqrt{\omega})$ for $n = 0, 1, \ldots$ . To compute $\mF_{\omega_1}^{(1,1)}$ in \eqref{oneeigsymm2}, we first obtain $R_\omega(\bx;\bx_0)$ from the definition in \eqref{Gkloc} along with the small argument asymptotics of $K_0(z)$


\begin{multline} \label{Gunitregular}
	R_\omega(\rho,\theta;\rho_0,\theta_0) = \frac{1}{2\pi}\left[-\gamma - \log\frac{\sqrt{\omega}}{2}   + \frac{\omega}{4}\left( -\log|\bx-\bx_0| + 1 - \gamma -\log\frac{\sqrt{\omega}}{2}  \right)|\bx-\bx_0|^2  \right] \\ - \frac{1}{2\pi}A_0I_0(\sqrt{\omega}\rho) - \frac{1}{\pi}\sum_{n  =1}^\infty \cos(n(\theta-\theta_0))A_n I_n(\sqrt{\omega}\rho) \,. 
\end{multline}

\noI In \eqref{Gunitregular}, $\gamma$ is Euler's constant. We next use \eqref{Gunitregular} to compute $F_{\omega_1}$, the gradient of the regular part with respect to $\bx$ evaluated at $\bx = \bx_0$. We have

\begin{multline} \label{Gunitregulargrad}
\nabla_\bx R_\omega(\bx;\bx_0) = \frac{\omega}{4\pi}\left(-\log|\bx-\bx_0| + \frac{1}{2} - \gamma - \log\frac{\sqrt{\omega}}{2} \right)(\bx-\bx_0) \\ - \frac{\sqrt{\omega}}{2\pi}A_0I_0^\prime(\sqrt{\omega}\rho)\be_\theta - \frac{\sqrt{\omega}}{\pi}\sum_{n  =1}^\infty \cos(n(\theta-\theta_0))A_n I_n^\prime(\sqrt{\omega}\rho)\be_\theta \,,
\end{multline}

where $\be_\theta \equiv (\cos\theta, \sin\theta)^T$. Setting $\bx = \bx_0$ in \eqref{Gunitregulargrad}, we obtain

\begin{equation} \label{Gunitregulargradself}
\bF_{\omega_1} \equiv \nabla_\bx R_\omega(\bx;\bx_0)\mid_{\bx=\bx_0} = - \frac{\sqrt{\omega}}{\pi}\sum_{n  =0}^\infty c_n A_n I_n^\prime(\sqrt{\omega}\rho_0)\be_{\theta_0} \,,
\end{equation}

where $c_n = 1/2$ ($c_n = 1$) when $n = 0$ ($n > 0)$. Observe in \eqref{Gunitregulargradself} that setting $\bx_0 = \bzero$ results in the gradient being $\bzero$, since $\bx_0 = \bzero$ is the equilibrium location of the spot.

Finally, to compute $\mF_{\omega_1}^{(1,1)}$ in \eqref{oneeigsymm2}, we use \eqref{gradGKdef} take the gradient of the first component of \eqref{Gunitregulargradself} with respect to $\bx_0$ using $\nabla_{\bx_0} = \be_{\theta_0}\partial_{\rho_0} + \rho_0^{-1}\be_{\theta_0}^\prime\partial_{\theta_0}$ to obtain 

\begin{multline} \label{Fomega1}
\nabla_{\bx_0} F_{\omega_1}^{(1)} =  -\be_{\theta_0} \cos\theta_0 \frac{\omega}{\pi} \sum_{n = 0}^\infty c_n \frac{K_n^\prime(\sqrt{\omega})}{I_n^\prime(\sqrt{\omega})}\left[\left( I_n^\prime(\sqrt{\omega}\rho_0) \right)^2 + I_n(\sqrt{\omega}\rho)I_n^{\prime\prime}(\sqrt{\omega}\rho_0) \right] \\ + \frac{1}{\rho_0}\be_{\theta_0}^\prime\sin\theta_0\frac{\sqrt{\omega}}{\pi}\sum_{n =0}^\infty c_n \frac{K_n^\prime(\sqrt{\omega})}{I_n^\prime(\sqrt{\omega})}I_n(\sqrt{\omega}\rho_0)I_n^\prime(\sqrt{\omega}\rho_0) \,,
\end{multline}

where $\be_{\theta_0}^\prime = (-\sin\theta, \cos\theta)$. Taking the limit as $\rho_0 \to 0^+$ in \eqref{Fomega1} yields

\begin{equation} \label{Fomega1_2}
	\left(\begin{array}{c} \mF_{\omega_1}^{(1,1)} \\ \mF_{\omega_1}^{(1,2)}  \end{array}\right) = \lim_{\rho_0 \to 0^+} \nabla_{\bx_0} F_{\omega_1}^{(1)} = -\frac{\omega}{4\pi}\left\lbrack \frac{K_0^\prime(\sqrt{\omega})}{I_0^\prime(\sqrt{\omega})} + \frac{K_1^\prime(\sqrt{\omega})}{I_1^\prime(\sqrt{\omega})} \right\rbrack \left(\begin{array}{c} 1 \\ 0 \end{array}\right) \,.
\end{equation}

\noI As expected, the second component of \eqref{Fomega1_2} is zero since the matrix $\mF_{\omega_1}$ must be diagonal due to the symmetry of the disk. It remains now to compute the Hessian term of the Neumann Green's function, $H_{11}^{(1,1)}$ and that of the Helmholtz Green's function, $H_{\omega_1}^{(1,1)}$.

The Neumann Green's function $G(\rho)$ satisfying \eqref{GNall} with source at the origin is given in polar coordinates by 

\BE \label{GNunit}
	G(\rho) = -\frac{1}{2\pi}\log\rho + \frac{\rho^2}{4\pi} - \frac{3}{8\pi} \,.
\EE

From \eqref{GNunit} and \eqref{GNloc}, we obtain that $H_{11}^{(1,1)} = (2\pi)^{-1}$. The Helmholtz Green's function $G_\omega(\rho)$ satisfying \eqref{Gk} with source at the origin is given in polar coordinates by

\BE \label{GHelmunit}
	G_\omega(\rho) = \frac{1}{2\pi}\left[K_0(\sqrt{\omega}\rho) - \frac{K_0^\prime(\sqrt{\omega})}{I_0^\prime(\sqrt{\omega})}I_0(\sqrt{\omega}\rho) \right] \,.
\EE

Using the small argument asymptotics of $K_0(z)$ and $I_0(z)$ in \eqref{GHelmunit}, we obtain from \eqref{Gkloc} the Hessian term

\BE \label{GHelmunithess}
	H_{\omega_1}^{(1,1)} = \frac{1}{\pi}\left[ -\frac{\omega}{4}\log\frac{\sqrt{\omega}}{2} + \frac{1-\gamma}{4}\omega \right] - \frac{\omega}{4\pi}\frac{K_0^\prime(\sqrt{\omega})}{I_0^\prime(\sqrt{\omega})} \,.
\EE

With $H_{\omega_1}^{(1,1)} = (2\pi)^{-1}$ and using \eqref{GHelmunithess} for $H_{\omega_1}^{(1,1)}$ and \eqref{Fomega1_2} for $\mF_{\omega_1}^{(1,1)}$, we obtain from \eqref{oneeigsymm2}

\BE \label{eigunit}
	2 + \omega\left[\log\left(\frac{e^\gamma \eps\sqrt{\omega}}{2}\right) - \frac{K_1^\prime(\sqrt{\omega})}{I_1^\prime(\sqrt{\omega})}\right] = \frac{\lambda_0k_{11} - \omega k_{21}}{S} \,.
\EE 

Setting $\lambda_0 = i\lambda_I$ and replacing $\omega$ with $i\lambda_I\tau_0$ in \eqref{eigunit}, we recover the complex equation given in (4.51) of \cite{xie2017moving} for the Hopf stability threshold $\tau_0$ and corresponding bifurcation frequency $\lambda_I$.

\subsection{The perturbed disk} \label{perturbeddisk}

In this section, we compute the leading order correction to the Neumann and Helmholtz Green's functions when $\Omega$ is the perturbed disk with boundary parametrized by $(x_1, x_2) = (1 + \sigma f(\theta))(\cos\theta,\sin\theta)$ with $\eps \ll \sigma \ll 1$ and some function $f(\theta)$ periodic over the interval $[0, 2\pi)$. We use this calculation to obtain analytic insight into how domain geometry impacts the preferred direction of oscillation at the onset of instability. Since $f$ can be expressed as a Fourier series, and the leading order effects of the perturbation are linear, it suffices to perform the calculation for individual modes $f(\theta) = \cos n\theta$ with $n = \mathbb{N}$. Analysis of perturbations of $\sin n\theta$ can can be recovered by replacing $\theta \to \theta - \pi/(2n)$.

We show in the following analysis that the $n=2$ mode impacts the bifurcation threshold, oscillation frequency, as well as oscillation mode at $\mO(\sigma)$. We then use this analysis to briefly show that the $n \neq 2$ effects enter only in higher orders $\sigma$.  We thus focus our calculations on the $n=2$ case for which the perturbed disk $\Omega$ is slightly elliptical in shape with the major axis aligned along the $x_1$ direction. For $n=2$, we show that the mode of oscillation along the $x_1$-axis is the first to become unstable as $\tau$ is increased. That is, we show that the threshold corresponding to the $(1,0)$ mode of oscillation is smaller than that for the unit disk, and that the threshold of the $(0,1)$ mode is larger than that of the unit disk.

The effect of the boundary perturbation on the localized variable $v(\bx)$ is exponentially small. As such, we need only expand 

\BE \label{usigma}
u \sim u_0(\rho) + \sigma u_1(\rho, \theta) \,,
\EE

and compute the boundary conditions for $u_1$ in terms of $u_0$. Proceeding, we find that an outward pointing normal vector on $\partial\Omega$ is

\BE \label{normvec}
	\mathbf{n} = (1+\sigma f(\theta)\be_\theta - \sigma f^\prime(\theta)\be_\theta^\prime \,.
\EE

In polar coordinates, the homogeneous boundary condition $\nabla u \cdot \mathbf{n} = 0$ becomes

\BE \label{pertbc}
	u_\rho - \frac{\sigma f^\prime}{(1+\sigma f)^2} u_\theta = 0 \,.
\EE

Substituting \eqref{usigma} into \eqref{pertbc} and expanding to first order in $\sigma$, we obtain the boundary conditions

\BE \label{pertbc2}
	u_{0\rho}(1) = 0 \,, \qquad u_{1\rho}(1,\theta) = - f(\theta) u_{0\rho\rho}(1) \,.
\EE

We begin with the expansion of the Neumann Green's function $G \sim G_0 + \sigma G_1$, where $G_0$ is the solution given in \eqref{GNunit} for the unperturbed unit disk $\Omega_0$ with source at the origin. From \eqref{pertbc2}, the boundary value problem for $G_1$ is then

\BE \label{G1bvp}
	\Delta G_{1\rho\rho} = 0 \,, \qquad G_{1\rho}(1,\theta) = -f(\theta)\frac{1}{\pi} \,; \qquad \int_{\Omega_0} \! G_1 \, d\Omega_0 = 0 \,,
\EE

with regularity at $\rho = 0$. When $f(\theta) = \cos  2\theta$, the solution to \eqref{G1bvp} is

\BE \label{G1sol}
	G_{1}(\rho,\theta) = -\frac{1}{2\pi}\rho^2\cos 2\theta \,,
\EE

or $G_1(x_1,x_2) = (2\pi)^{-1}(x_2^2-x_1^2)$ in Cartesian coordinates. To leading order in $\sigma$, we therefore have that the Hessian matrix $H_{11}$ in \eqref{oneeigdef} is given by

\BE \label{H11sigma}
 H_{11} \sim	\frac{1}{2\pi}\mI_2 + \sigma \frac{1}{\pi}\left(\begin{array}{cc} -1 & 0 \\ 0 & 1\end{array}\right) \,,
\EE

where the leading order term in \eqref{H11sigma} is computed from Neumann's Green's function of the unperturbed disk given in \eqref{GNunit}.

For the Helmholtz Green's function $G_\omega$ of \eqref{Gkeq}, we expand $G_{\omega} = G_{\omega_0} + \sigma G_{\omega_1}$, where $G_{\omega_0}$ is the leading order solution when $\rho > \rho_0$ given by \cite{chen2011stability}

\BE \label{Gomega0}
	G_{\omega_0}(\rho,\theta; \rho_0, \theta_0) = \frac{1}{2\pi} \sum_{n = -1}^\infty e^{-in(\theta-\theta_0)} F_n(\rho) I_n(\sqrt{\omega}\rho_0) \,,
\EE

where $F_n(\rho) \equiv K_n(\sqrt{\omega}\rho) - K_n^\prime(\sqrt{\omega})I_n(\sqrt{\omega}\rho)/I_n^\prime(\sqrt{\omega})$. We note that \eqref{Gomega0} is equivalent to \eqref{Gunit} when $\rho > \rho_0$. While \eqref{Gunit} has the singularity extracted analytically from the sum, \eqref{Gomega0} is more convenient to work with when $\bx$ is not near $\bx_0$. From \eqref{pertbc2}, the boundary value problem for $G_1$ is then

\BE \label{Gw1bvp}
		\Delta G_{\omega_1 \rho\rho} - G_{\omega_1} = 0 \,, \qquad G_{\omega_1}(1,\theta) = \frac{1}{2\pi} \textrm{Re}\left[\sum_{n = -\infty}^\infty e^{-i(n-2)\theta + in\theta_0} F_n^{\prime\prime}(1) I_n(\sqrt{\omega}\rho_0)   \right] \,.
\EE

The solution to \eqref{Gw1bvp} is

\BE \label{Gw1sol}
	G_{\omega_1}(\rho,\theta;\rho_0,\theta_0) = \frac{1}{2\pi\sqrt{\omega}}\textrm{Re}\left[\sum_{n = -\infty}^\infty e^{-i(n-2)\theta + in\theta_0} F_n^{\prime\prime}(1) I_n(\sqrt{\omega}\rho_0)\frac{I_{n-2}(\sqrt{\omega}\rho)}{I_{n-2}^\prime(\sqrt{\omega})}  \right] \,.
\EE

Taking only the real part for $G_{\omega_1}$, we obtain

\BE \label{Gw1sol2}
	G_{\omega_1}(\rho,\theta;\rho_0,\theta_0) = \frac{1}{2\pi\sqrt{\omega}}\left[\sum_{n = -\infty}^\infty \cos\left[(n-2)\theta - n\theta_0\right] F_n^{\prime\prime}(1) I_n(\sqrt{\omega}\rho_0)\frac{I_{n-2}(\sqrt{\omega}\rho)}{I_{n-2}^\prime(\sqrt{\omega})}  \right] \,.
\EE

We now compute the contributions of the domain perturbation to the quantities  $\mF_{\omega_1}^{(k,k)}$ and $H_{\omega_1}^{(k,k)}$, $k = 1, 2$,  of \eqref{oneeigsymm1}. We first compute the correction to $\bF_{\omega}$, the gradient of the regular part evaluated at $\bx_0$. Since $G_{\omega_1}$ contains no singularities, we use $\nabla_{\bx} = \be_{\theta}\partial_{\rho} + \rho^{-1}\be_{\theta}^\prime\partial_{\theta}$ to obtain
	
\bes
\BE \label{gradGw1}	
	\bg \equiv \nabla_{\bx}G_{\omega_1} \mid_{\bx=\bx_0} = \be_{\theta_0} \cos 2\theta_0 A(\rho_0) + \be_{\theta_0}^\prime \sin 2\theta_0 B(\rho_0) \,,
\EE

where

\BE \label{gradGw1def}	
\begin{gathered}
	A(\rho_0) \equiv \frac{1}{2\pi} \sum_{n=-\infty}^\infty \frac{F_n^{\prime\prime}(1)}{I_{n-2}^\prime(\sqrt{\omega})}I_n(\sqrt{\omega}\rho_0)I_{n-2}^\prime(\sqrt{\omega}\rho_0) \,, \\
B(\rho_0) \equiv \frac{1}{2 \pi \sqrt{\omega} \rho_0} \sum_{n=-\infty}^\infty (n-2) \frac{F_n^{\prime\prime}(1)}{I_{n-2}^\prime(\sqrt{\omega})}I_n(\sqrt{\omega}\rho_0)I_{n-2}(\sqrt{\omega}\rho_0) \,.
\end{gathered}
\EE
\ees

The leading order correction to $\mF_\omega$ is given by 
$\lim_{\rho \to 0^+} \left(\begin{array}{cc} \nabla_{\bx_0} \bg^{(1)} & \nabla_{\bx_0} \bg^{(2)} \end{array}\right)$, where $\bg^{(1)}$ and $\bg^{(2)}$ are the first and second components of the vector $\bg$ defined in \eqref{gradGw1}. We next use the fact that for $|z| \ll 1$, we have

\BE \label{limbes}
\begin{gathered}
	I_n(z)I_{n-2}^\prime(z) \sim \begin{cases} 
      \frac{z}{4} & n = 0, 1 \\
      \mo(z^2) & \textrm{else} 
   \end{cases} \,; \\
	I_n(z)I_{n-2}(z) \sim \begin{cases} 
      \frac{z^2}{8} & n = 0, 2 \\
			\frac{z^2}{4} & n = 1 \\
      \mo(z^3) & \textrm{else} 
   \end{cases} \,.
\end{gathered}
\EE

We thus obtain that

\BE \label{Fcoor2}
	\lim_{\rho \to 0^+} \left(\begin{array}{cc} \nabla_{\bx_0} \bg^{(1)}  &\nabla_{\bx_0} \bg^{(2)} \end{array}\right) = \frac{1}{8\pi}\sqrt{\omega}\left[\frac{F_0^{\prime\prime}(1)}{I_{-2}^\prime(\sqrt{\omega})} + \frac{F_1^{\prime\prime}(1)}{I_{-1}^\prime(\sqrt{\omega})}\right]\left(\begin{array}{cc} 1 & 0 \\ 0 & -1\end{array}\right) \,. 
\EE

To leading order in $\sigma$, the matrix $\mF_{\omega_1}$ in \eqref{oneeigdef} is given by

\BE \label{Fwsigma}
	\mF_{\omega_1} \sim -\frac{\omega}{4\pi}\left\lbrack \frac{K_0^\prime(\sqrt{\omega})}{I_0^\prime(\sqrt{\omega})} + \frac{K_1^\prime(\sqrt{\omega})}{I_1^\prime(\sqrt{\omega})} \right\rbrack \mI_2 + \sigma \frac{1}{8\pi}\sqrt{\omega}\left[\frac{F_0^{\prime\prime}(1)}{I_{-2}(\sqrt{\omega})} + \frac{F_1^{\prime\prime}(1)}{I_{-1}(\sqrt{\omega})}\right]\left(\begin{array}{cc} 1 & 0 \\ 0 & -1\end{array}\right) \,,
\EE

where the leading order term of \eqref{Fwsigma} was computed in \eqref{Fomega1_2}.

Finally, for the correction to the Hessian term $H_{\omega_1}$ of \eqref{oneeigdef}, we set $(\rho,\theta) = (\rho_0,\theta_0)$ in \eqref{Gw1sol2} and let $\rho \to 0^{+}$ to obtain

\BE \label{Gw1sol2hess1}
		G_{\omega_1}(\bx_0;\bx_0) \sim \frac{\sqrt{\omega}}{8\pi}\left[\frac{F_0^{\prime\prime}(1)}{I_{-2}^\prime(\sqrt{\omega})} + \frac{F_1^{\prime\prime}(1)}{I_{-1}^\prime(\sqrt{\omega})}   \right] (x_2^2 - x_1^2) \,
\EE

yielding the two-term expansion for $H_{\omega_1}$

\begin{multline} \label{Hw1sigma}
	H_{\omega_1} \sim \left\lbrace \frac{1}{\pi}\left[ -\frac{\omega}{4}\log\frac{\sqrt{\omega}}{2} + \frac{1-\gamma}{4}\omega \right] - \frac{\omega}{4\pi}\frac{K_0^\prime(\sqrt{\omega})}{I_0^\prime(\sqrt{\omega})}\right\rbrace \mI_2 + \\  \sigma\frac{\sqrt{\omega}}{4\pi}\left[\frac{F_0^{\prime\prime}(1)}{I_{-2}^\prime(\sqrt{\omega})} + \frac{F_1^{\prime\prime}(1)}{I_{-1}^\prime(\sqrt{\omega})} \right]\left(\begin{array}{cc} -1 & 0 \\ 0 & 1\end{array}\right) \,,
\end{multline}

where the leading order term in \eqref{Hw1sigma} corresponds to that for the unperturbed unit disk computed in \eqref{GHelmunithess}. 

We now substitute \eqref{H11sigma} for $H_{11}$, \eqref{Fwsigma} for $\mF_{\omega_1}$, and \eqref{Hw1sigma} for $H_{\omega_1}$ into \eqref{oneeigsymm1}. We perturb the eigenvalue $\lambda_0 = i(\lambda_{I_0} + \sigma \lambda_{I_1})$ and the stability threshold $\tau_0 = \hat{\tau}_0 + \sigma\hat{\tau}_1$ so that $\omega = \lambda_0\tau_0 \equiv \omega_0 + \sigma \omega_1$, where $\omega_0 = i\lambda_{I_0}\hat{\tau}_0$ and $\omega_{1} = i(\lambda_{I_1}\hat{\tau}_0 + \lambda_{I_0}\hat{\tau}_1)$. The calculation below will show that $\hat{\tau}_1 < 0$ for the $(1,0)$ oscillation mode while $\hat{\tau}_1 > 0$ for the $(0,1)$ oscillation mode, implying that the preferred direction of oscillation is along the major axis of the perturbed unit disk.

From \eqref{oneeigsymm1}, we obtain the leading order eigenvalue problem \eqref{eigunit} for $\omega_0$ and $\lambda_{I_0}$ corresponding to the unperturbed problem, while at $\mO(\sigma)$, we obtain the two uncouple equations for $\omega_1$ and $\lambda_{I_0}$

\bes \label{Osigmaeig}
\begin{multline} \label{Osigmaeig1}
	\omega_1\left[\log\frac{e^\gamma\eps}{2} + \frac{1}{2}\log\omega_0 - \frac{K_1^\prime(\sqrt{\omega_0})}{I_1^\prime(\sqrt{\omega_0})} \right] + \omega_0\left[\frac{1}{2}\frac{\omega_1}{\omega_0} - \omega_1 Q(\sqrt{\omega_0}) \right]   \\  - \frac{1}{\pi} + \frac{3}{8\pi}\sqrt{\omega_0} \left[ \frac{F_0^{\prime\prime}(1)}{I_{-2}^\prime(\sqrt{\omega_0})} + \frac{F_1^{\prime\prime}(1)}{I_{-1}^\prime(\sqrt{\omega_0})}  \right] = \frac{ik_{11}\lambda_{I_1} - k_{21}\omega_1}{S} \,,
\end{multline}
\begin{multline} \label{Osigmaeig2}
	\omega_1\left[\log\frac{e^\gamma\eps}{2} + \frac{1}{2}\log\omega_0 - \frac{K_1^\prime(\sqrt{\omega_0})}{I_1^\prime(\sqrt{\omega_0})} \right] + \omega_0\left[\frac{1}{2}\frac{\omega_1}{\omega_0} - \omega_1 Q(\sqrt{\omega_0}) \right]   \\ + \frac{1}{\pi} - \frac{3}{8\pi}\sqrt{\omega_0} \left[ \frac{F_0^{\prime\prime}(1)}{I_{-2}^\prime(\sqrt{\omega_0})} + \frac{F_1^{\prime\prime}(1)}{I_{-1}^\prime(\sqrt{\omega_0})}  \right] = \frac{ik_{11}\lambda_{I_1} - k_{21}\omega_1 }{S} \,,
\end{multline}

where the function $Q(z)$ is defined by

\BE \label{Osigmaeig1def}
	Q(z) = \frac{1}{2z}\frac{(z^2+1)\left[K_1(z)I_0(z) + K_0(z)I_1(z)\right]}{(zI_0(z) - I_1(z))^2} \,.
\EE
\ees

The signs in \eqref{Osigmaeig1} correspond to that of the $(1,0)$ mode of oscillation, while the signs of \eqref{Osigmaeig2} correspond to that of  the $(0,1)$ mode. Rearranging \eqref{Osigmaeig1} and \eqref{Osigmaeig2}

\bes \label{Osigmaeigall}
\begin{multline} \label{Osigmaeig3}
	\omega_1\left[\log\frac{e^\gamma\eps}{2} + \frac{1}{2}\log\omega_0 - \frac{K_1^\prime(\sqrt{\omega_0})}{I_1^\prime(\sqrt{\omega_0})} + \frac{1}{2} - \omega_0Q(\sqrt{\omega_0}) + \frac{k_{21}}{S}\right] - \frac{ik_{11}\lambda_{I_1}}{S} \\ = \frac{1}{\pi} -  \frac{3}{8\pi}\sqrt{\omega_0} \left[ \frac{F_0^{\prime\prime}(1)}{I_{-2}^\prime(\sqrt{\omega_0})} + \frac{F_1^{\prime\prime}(1)}{I_{-1}^\prime(\sqrt{\omega_0})}  \right]  \,,
\end{multline}
\begin{multline} \label{Osigmaeig4}
	\omega_1\left[\log\frac{e^\gamma\eps}{2} + \frac{1}{2}\log\omega_0 - \frac{K_1^\prime(\sqrt{\omega_0})}{I_1^\prime(\sqrt{\omega_0})} + \frac{1}{2} - \omega_0Q(\sqrt{\omega_0}) + \frac{k_{21}}{S}\right] - \frac{ik_{11}\lambda_{I_1}}{S} \\ = -\frac{1}{\pi} +  \frac{3}{8\pi}\sqrt{\omega_0} \left[ \frac{F_0^{\prime\prime}(1)}{I_{-2}^\prime(\sqrt{\omega_0})} + \frac{F_1^{\prime\prime}(1)}{I_{-1}^\prime(\sqrt{\omega_0})}  \right]  \,.
\end{multline}
\ees

Noting that $\omega_1 = i(\lambda_{I_1}\hat{\tau}_0 + \lambda_{I_0}\hat{\tau}_1)$ is pure imaginary, we  equate the real and imaginary parts of \eqref{Osigmaeigall} to obtain

\bes \label{Osigmaeigall2}
\begin{multline} \label{Osigmaeig5}
	-(\lambda_{I_1}\hat{\tau}_0 + \lambda_{I_0}\hat{\tau}_1)\textrm{Im}\left\lbrace\frac{1}{2}\log\omega_0 - \frac{K_1^\prime(\sqrt{\omega_0})}{I_1^\prime(\sqrt{\omega_0})} - \omega_0Q(\sqrt{\omega_0}) \right\rbrace  \\ = \pm \frac{1}{\pi} \mp  \frac{3}{8\pi} \textrm{Re}\left\lbrace \sqrt{\omega_0} \left[ \frac{F_0^{\prime\prime}(1)}{I_{-2}^\prime(\sqrt{\omega_0})} + \frac{F_1^{\prime\prime}(1)}{I_{-1}^\prime(\sqrt{\omega_0})}  \right] \right\rbrace  \,,
\end{multline}
\begin{multline} \label{Osigmaeig6}
	(\lambda_{I_1}\hat{\tau}_0 + \lambda_{I_0}\hat{\tau}_1)\textrm{Re}\left\lbrace \frac{k_{11}}{S\hat{\tau}_0} + \frac{1}{2} - \omega_0Q(\sqrt{\omega_0}) \right\rbrace - \frac{k_{11}\lambda_{I_1}}{S} \\ =  \mp  \frac{3}{8\pi} \textrm{Im}\left\lbrace\sqrt{\omega_0} \left[ \frac{F_0^{\prime\prime}(1)}{I_{-2}^\prime(\sqrt{\omega_0})} + \frac{F_1^{\prime\prime}(1)}{I_{-1}^\prime(\sqrt{\omega_0})}  \right]  \right\rbrace \,,
\end{multline}
\ees

where, in \eqref{Osigmaeig6}, we have used (4.52b) of \cite{xie2017moving} to obtain that the leading order threshold $\hat{\tau}_0$ satisfies the relationship

\BE \label{tau0hat}
	\textrm{Re}\left\lbrace \log\frac{e^\gamma \eps}{2} + \frac{1}{2}\log\omega_0 -\frac{K_1^\prime(\sqrt{\omega_0})}{I_1^\prime(\sqrt{\omega_0}) } + \frac{k_{21}}{S} \right\rbrace = \frac{k_{11}}{S\hat{\tau}_0}  \,.
\EE 

Observe that $\hat{\tau}_0 \sim 1/|\log\eps| \ll 1$. In \eqref{Osigmaeigall2}, the top (bottom) sign corresponds to the $(1,0)$  ($(0,1)$) oscillation mode. Solving for the quantity $\lambda_{I_1}\hat{\tau}_0 + \lambda_{I_0}\hat{\tau}_1$ in \eqref{Osigmaeig5}, we obtain

\BE \label{omega1}
   (\lambda_{I_1}\hat{\tau}_0 + \lambda_{I_0}\hat{\tau}_1) \equiv \tilde{\omega}_1 = \frac{\mp\frac{1}{\pi} \pm  \frac{3}{8\pi} \textrm{Re}\left\lbrace \sqrt{\omega_0} \left[ \frac{F_0^{\prime\prime}(1)}{I_{-2}^\prime(\sqrt{\omega_0})} + \frac{F_1^{\prime\prime}(1)}{I_{-1}^\prime(\sqrt{\omega_0})}  \right] \right\rbrace}{\textrm{Im}\left\lbrace\frac{1}{2}\log\omega_0 - \frac{K_1^\prime(\sqrt{\omega_0})}{I_1^\prime(\sqrt{\omega_0})} - \omega_0Q(\sqrt{\omega_0}) \right\rbrace} \,.
\EE

We note that there are no parameters on the right-hand side of \eqref{omega1}, as the numerical value for $\omega_0$ was given in (4.52c) of \cite{xie2017moving} as $\omega_0 = 3.02603687i$ and shown to be independent of parameters of the original PDE system \eqref{schnak}. From \eqref{Osigmaeig6}, we solve for $\lambda_{I_1}$ to obtain

\begin{equation} \label{lambda1}
	\lambda_{I_1} = \frac{\tilde{\omega}_1}{\hat{\tau}_0} + \frac{S}{k_{11}}\left[\tilde{\omega}_1\textrm{Re}\left\lbrace \frac{1}{2} - \omega_0Q(\sqrt{\omega_0}) \right\rbrace \pm  \frac{3}{8\pi} \textrm{Im}\left\lbrace\sqrt{\omega_0} \left[ \frac{F_0^{\prime\prime}(1)}{I_{-2}^\prime(\sqrt{\omega_0})} + \frac{F_1^{\prime\prime}(1)}{I_{-1}^\prime(\sqrt{\omega_0})}  \right]  \right\rbrace \right]  \,.
\end{equation}

Notice that as $\eps \to 0$, since $\hat{\tau}_0 \sim 1/|\log\eps|$, $\lambda_{I_1} \sim \tilde{\omega}_1|\log\eps|$. That is, $\lambda_{I_1}$ has the same sign as $\tilde{\omega}_1$ as $\eps \to 0$.

Finally, with $\hat{\tau}_1 = (\tilde{\omega}_1 - \lambda_{I_1}\hat{\tau}_0)/\lambda_{I_0}$, we use \eqref{lambda1} for $\lambda_{I_1}$ to obtain

\BE \label{tau1}
	\hat{\tau}_1 = -\frac{S\hat{\tau}_0}{k_{11}\lambda_{I_0}}\left[ \tilde{\omega_1}  \textrm{Re}\left\lbrace \frac{1}{2} - \omega_0Q(\sqrt{\omega_0}) \right\rbrace \pm  \frac{3}{8\pi} \textrm{Im}\left\lbrace\sqrt{\omega_0} \left[ \frac{F_0^{\prime\prime}(1)}{I_{-2}^\prime(\sqrt{\omega_0})} + \frac{F_1^{\prime\prime}(1)}{I_{-1}^\prime(\sqrt{\omega_0})}  \right]  \right\rbrace  \right] \,.
\EE

Noting that the spot strength $S$, the leading order threshold $\hat{\tau}_0$, and leading order frequency $\lambda_{I_0}$ are all positive while $k_{11} < 0$ from Fig. \ref{kappafig}, we use $\omega_0 = 3.02603687i$ and the upper signs in \eqref{omega1} and \eqref{tau1} to find that $\hat{\tau}_1 \approx -(1.083)\hat{\tau}_0/(|k_{11}|\lambda_{I_0}) < 0 $ while $\hat{\tau}_1 \approx (1.083)\hat{\tau}_0/(|k_{11}|\lambda_{I_0}) > 0 $. We thus conclude that the for the perturbed disk with radius $r = 1 + \sigma\cos 2\theta$ with $0 < \sigma \ll 1$, preferred mode of oscillation is along the major axis, with corresponding threshold $\hat{\tau} \sim \hat{\tau}_0 + \sigma\hat{\tau}_1$ less than that for the unit disk. The threshold corresponding to oscillation along the minor axis is larger than that for the unit disk. Furthermore, since $\tilde{\omega}_1 < 0$ for the $(1,0)$ mode, the frequency of oscillation is less than that for the unit disk. Likewise,  $\tilde{\omega}_1 > 0$ for the $(0,1)$ mode so that the oscillation frequency along that direction increases with the perturbation of the unit disk.

A perturbation of the form $f(\theta) = a_2\cos 2\theta + b_2\sin 2\theta$ can be written as $f(\theta) = C_2\cos (2\theta -\phi)$, where $C_2 = \sqrt{a_2^2 + b_2^2}$ and $\cos\phi = a_2/C_2$, $\sin\phi = b_2/C_2$. This simply constitutes a counterclockwise rotation of the perturbation by $\phi/2$. In this case, the direction of oscillation is along the line that makes an angle $\phi/2$ with the $x_1$ axis.

We now briefly comment on the effect of perturbations of the form $f(\theta) = \cos m\theta$ for $m \neq 2$. The $m = 1$ perturbation constitutes a translation of the unit disk by $\sigma$ in the $x_1$ direction along with an $\mO(\sigma^2)$ deformation away from circular geometry. That is, the parametric equation for this perturbed disk is $x_1 = \sigma + (1+ \sigma^2 g(\theta))\cos\theta$ and $x_2 = (1+ \sigma^2 g(\theta))\sin\theta$ with $g(\theta) = (1/4)(1-\cos 2\theta) + \mO(\sigma^2)$. The geometry deformation to leading order is thus effectively a mode-$2$ perturbation with negative coefficient. From the above analysis, the $m=1$ perturbation thus selects the $(0,1)$ oscillation mode as the dominant mode. We emphasize, however, that this is an $\mO(\sigma^2)$ effect in contrast to the $\mO(\sigma)$ effect of the $m = 2$ mode. A similar argument asserts that the effect of $f(\theta) = \sin\theta$ is also $\mO(\sigma^2)$.

When $m \geq 3$, the leading order correction $G_1(\rho,\theta)$ to the Neumann Green's function in \eqref{G1sol} would be replaced by an $\mO(\rho^m)$ function, which has a zero Hessian at the origin. It therefore does not contribute to the eigenvalue problem \eqref{oneeigsymm2}. For the Helmholtz Green's function, we have that $I_m(z) \sim z^{|m|}$ for $|z| \ll 1$ so that $I_{n}(z)I_{n-m}(z) \sim z^{\alpha}$ where $\alpha = |n| + |n-m| \geq |m| $ while $I_{n}(z)I_{n-m}^\prime(z) \sim z^{\alpha}$ where $\alpha = |n| + |n-m-1| \geq |m-1|$. From \eqref{limbes}, we observe modes $m \geq 3$ do not contribute towards $\mF_{\omega_1}$ in \eqref{Fwsigma} and $H_{\omega_1}$ in \eqref{Hw1sigma}.

Since the leading order corrections to the Neumann and Helmholtz Green's functions are linear in the perturbation $f(\theta)$, we can conclude that for a $2\pi$-periodic function $f(\theta)$, the leading order contribution to the eigenvalue problem comes from the coefficients $a_2 \equiv \pi^{-1}\int_0^{2\pi} \!f(\theta)\cos 2\theta \, d\theta$ and $b_2 \equiv \pi^{-1}\int_0^{2\pi} \!f(\theta)\sin 2\theta \, d\theta$ in the Fourier series of $f$. In particular, if $a_2  = b_2= 0$, the bifurcation threshold (to leading order in $\sigma$) will remain unchanged from that of the unit disk. When $a_2$ and $b_2$ are not both zero, the bifurcation threshold will decrease, and the preferred direction of oscillation will be along a line that makes an angle $\phi/2$ with the $x_1$-axis, where $\cos\phi = a_2/\sqrt{a_2^2 + b_2^2}$, $\sin\phi = b_2/\sqrt{a_2^2+b_2^2}$. Furthermore, the corresponding bifurcation frequency will also decrease.

While we have not performed this analysis for near-square rectangle with edge lengths $L$ and $L+\sigma$ for $\eps \ll \sigma \ll 1$, our analysis for the perturbed disk predicts that the preferred oscillation direction in a rectangle would be in the direction parallel to the longer edge. We show a numerical example of this scenario, along with several others illustrating our single- and multi-spot theory, in the next section.

%

\section{Numerical validation} \label{numerics}
 In this section, we numerically validate our theoretical result \eqref{finaleigblock} by solving the full time-dependent PDE system \eqref{schnak} using the finite element solver FlexPDE 7 \cite{flex}. For single- and multi-spot equilibria of \eqref{schnak} on various domain geometries, we verify both the Hopf bifurcation threshold $\tau = \eps^{-2}\hat{\tau}^*$ for the onset of oscillatory instabilities as well as the direction(s) of oscillations, where $\hat{\tau}^*$ is the minimum of all $\hat{\tau}$ satisfying the complex eigenvalue problem \eqref{finaleigblock}.

Before we describe our results, we first outline our procedures. Initial $N$-spot equilibrium states for which we test stability are obtained by initializing an $N$-bump pattern in \eqref{schnak} with $\tau$ set well below any Hopf thresholds. We then evolve \eqref{schnak} until the time $t$ is sufficiently large that changes in solution are no longer observed; observe that steady state solutions of \eqref{schnak} are unaffected by the value $\tau$. Using this equilibrium solution as an initial condition, we trial various of values of $\tau$ to test stability; we say that the numerical (or ``exact'') value of the Hopf bifurcation threshold is $\hat{\tau}_f$ if no oscillations are observed when $\red{\tau_0} < \hat{\tau}_f-0.005$ and oscillations are observed when $\red{\tau_0} > \hat{\tau}_f + 0.005$. 

To compute $\hat{\tau}$ from the matrix-eigenvalue problem, \eqref{det0}, we require quantities associated with the Neumann and Helmholtz Green's function of \eqref{GNall} and \eqref{Gk}, respectively, for the domains $\Omega$ that we consider. For the unit disk and rectangle, analytic formulas for both are given in \cite{chen2011stability}. For the half disk, we simply employ the unit disk formula with the method of images to obtain the reflective boundary condition on the straight segment of the half-disk. For the more complex geometries, we employ the finite element solver from MATLAB's PDE Toolbox. In the implementation, we solve a regular equation for $G - G^{\textrm{free}}$ with nonhomogeneous Neumann conditions, where $G$ is the desired Green's function in $\Omega$ and $G^{\textrm{free}}$ is the corresponding free space Green's function. Finally, the equilibrium locations $\bx_j$ and corresponding spot strengths $S_j$ are obtained by simultaneously solving the $3N+1$ system of nonlinear equations \eqref{Sj}, \eqref{eqbm1}, and \eqref{matchj}.

\red{We remark that there are as many as $2N$ distinct pairs of solutions $(\hat{\omega}, \hat{\lambda_I})$ to the nonlinear system \eqref{det0}. In the $N$-spot numerical experiments below, we were were not always able to find all thresholds when $N$ was sufficiently large. We make a note of this below where relevant. In all cases, however, the smallest threshold we found corresponded to the most unstable mode as observed numerically in the PDE simulations.}

\subsection{Hopf bifurcation of a single spot}
 In this subsection, we investigate the Hopf bifurcation of small eigenvalues of a single spot solution to \eqref{schnak} in three different domains $\Omega$. The experiments are:
 \begin{enumerate}
		\item $\Omega$ is a perturbed unit disk with radius 
		\begin{equation} \label{pertf}
		\begin{gathered}
		r = 1 + \sigma f(\theta) \,; \qquad \theta \in \lbrack 0, 2\pi)  \,, \quad \sigma  = 20\eps \,; \\
		f(\theta) = \frac{1}{4}\left \lbrack 0.1\cos\theta + a_2\cos 2\theta + \cos 3\theta + \cos 4\theta \right\rbrack \,.
		\end{gathered}
		\end{equation}
		In one numerical simulation, we take $a_2 = 0.1$, while we take $a_2 = -0.1$ in the other simulation. We confirm our prediction of \S \ref{perturbeddisk} for the impact of the $a_2$ coefficient on the dominant mode of oscillation.
    \item  $\Omega$ is a half-disk with radius $1$ (Fig.~\ref{halfdisk1}). The predicted Hopf threshold is $\hat{\tau}^* = 0.0712$ while the numerical result was found to be $\hat{\tau}_f = 0.072$.  
    \item  $\Omega$ is a rectangle of width $2$ and length $1$ (see Fig.~\ref{rectangle1}). The predicted Hopf threshold is $\hat{\tau}^* = 0.0615$, while the numerical threshold was found to be $\hat{\tau}_f = 0.062$.  
    \item $\Omega$ is an asymmetric, non-simply connected domain consisting of the same rectangle as that in Fig.~\ref{rectangle1} with two differently sized holes in the shape of disks (see Fig.~\ref{rectangle2}).  All boundary conditions are reflective. The predicted Hopf threshold is $\hat{\tau} = 0.065$ while the numerical threshold was found to be $\hat{\tau}_f = 0.0615$.
 \end{enumerate}

In each of the figures below showing snapshots of the solution for the localized activator component at various points in time, the red arrows indicate the direction of motion of the spot(s) at the given instant. We note that since the strength of the spot is a function of its location, spot-splitting may occur during the course of oscillations. In the multi-spot solutions, spot annihilation events may occur when spot distances decrease due to oscillation. We exclude illustrations of these phenomena and focus only on the oscillatory dynamics.

\begin{figure}[!htb]
    \centering
    \includegraphics[width=0.98\textwidth]{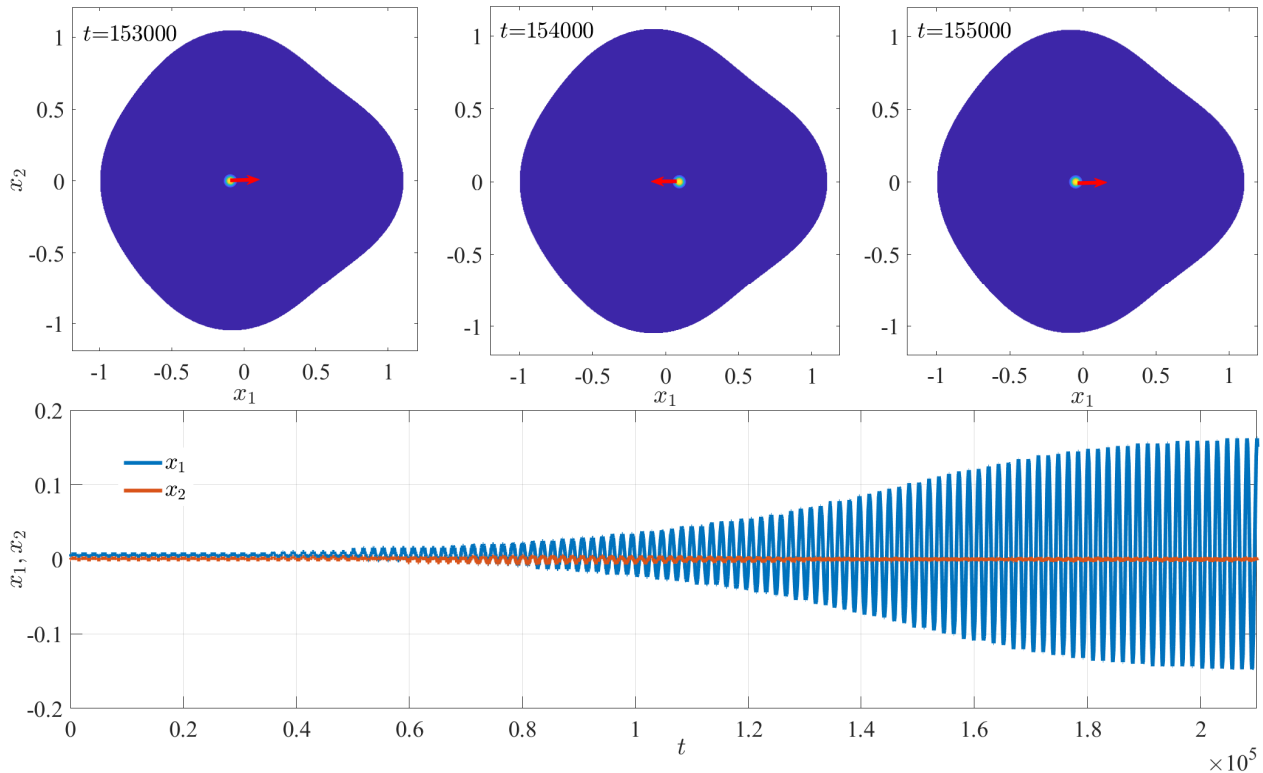}		
       \caption{Experiment 1 -- Hopf bifurcation of a single spot in a perturbed disk of the form \eqref{pertf} with $a_2 = 0.1 > 0$. In the top half, we indicate with arrows the direction of spot motion at the particular time. In the bottom half, we plot the $x_1$ and $x_2$ coordinates of the spot location as functions of time. For increasing time, the horizontal oscillation mode emerges as dominant. The other parameters are $S = 4$, $\eps = 0.01$, and $\tau_0 = 0.084$. }
    \label{a2positive}
\end{figure}

\textbf{Experiment 1.} \, In this experiment, we verify the calculations of \S \ref{perturbeddisk} on how the perturbation of the disk impacts the dominant mode of oscillation. For the perturbed disk of \eqref{pertf} with $a_2 = 0.1 > 0$, Fig. \ref{a2positive} shows the oscillations of the spot at onset as well as in long-time. In particular, we observe that the $(1,0)$ oscillation mode (i.e., the horizontal mode) is the dominant mode, consistent with the result of \S \ref{perturbeddisk}. 

\begin{figure}[!htb]
    \centering
    \includegraphics[width=0.98\textwidth]{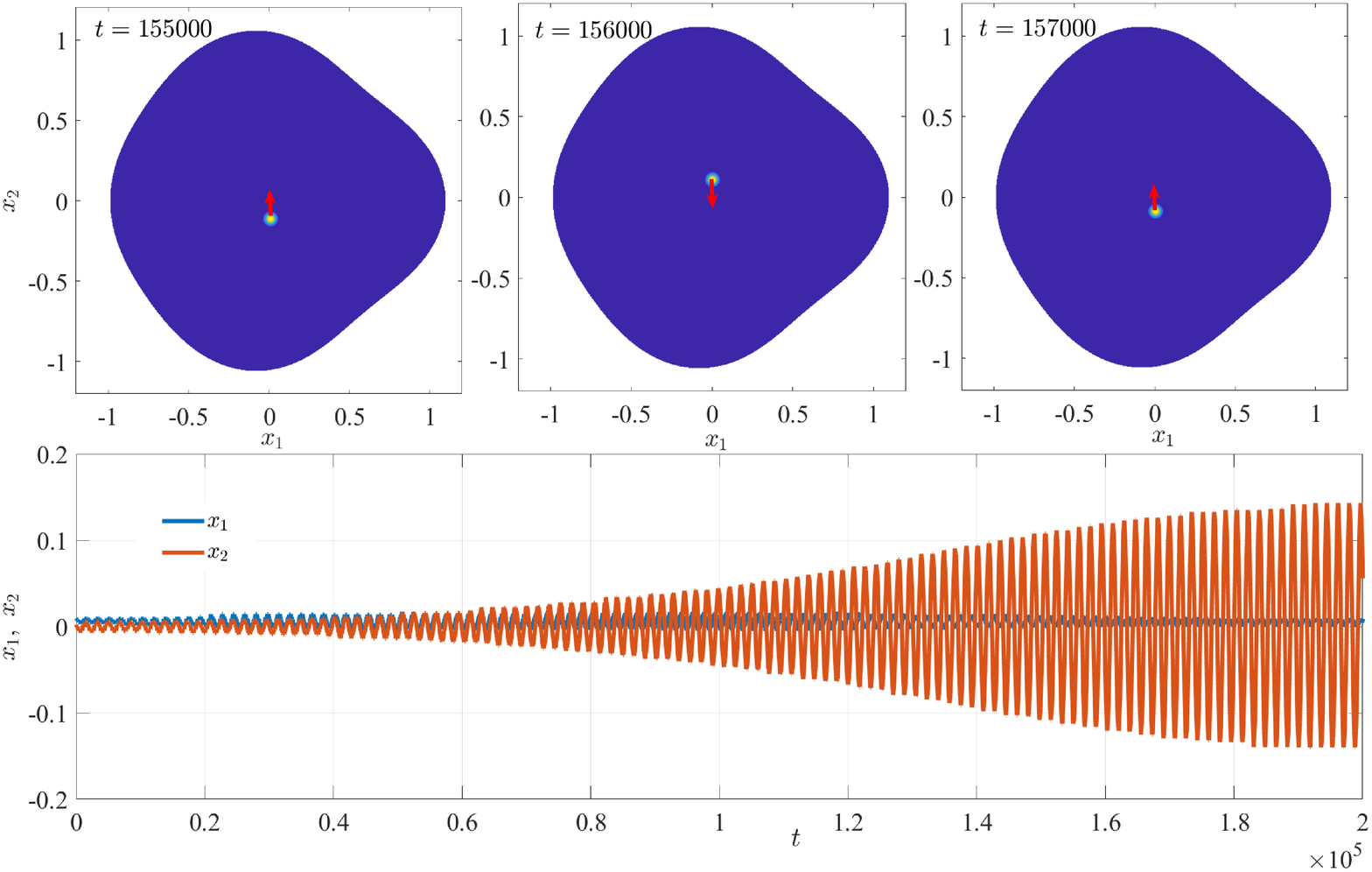}		
       \caption{Experiment 1 -- Hopf bifurcation of a single spot in a perturbed disk of the form \eqref{pertf} with $a_2 = -0.1 < 0$. In the top half, we indicate with arrows the direction of spot motion at the particular time. In the bottom half, we plot the $x_1$ and $x_2$ coordinates of the spot location as functions of time. For increasing time, the vertical oscillation mode emerges as dominant. The other parameters are $S = 4$, $\eps = 0.01$, and $\tau_0 = 0.084$. }
    \label{a2negative}
\end{figure}

\begin{figure}[!htb]
    \centering
    \includegraphics[width=0.98\textwidth]{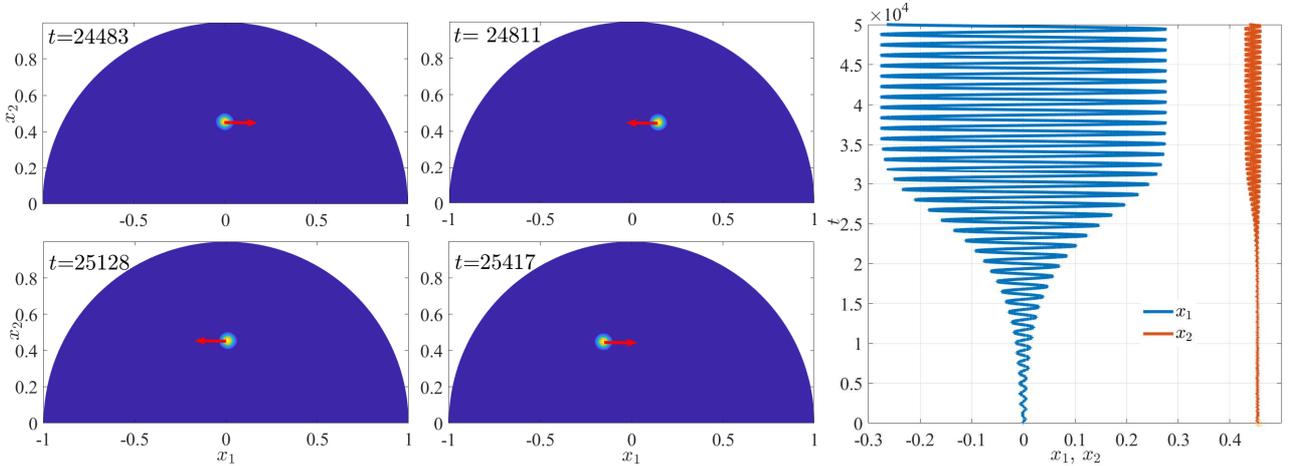}		
       \caption{Experiment 2 -- numerical simulations performed for $\tau_0=0.0725,~\varepsilon=0.01$ and $S= 4$. On the left, we show oscillatory motion of the spot center as $\tau_0$ exceeds the Hopf bifurcation threshold of $\hat{\tau}^* = 0.072$. Red arrows indicate the direction of motion. On the right, we plot coordinates of the center of the spot, where the blue (red) curve is the $x_1$-coordinate ($x_2$-coordinate). The saturation of the oscillation amplitudes indicates that the Hopf bifurcation is supercritical, with the initial horizontal oscillations leading to a stable orbit with nonzero $x_1$ and $x_2$ components.}
    \label{halfdisk1}
\end{figure}

In Fig. \ref{a2negative}, we set $a_2 = -0.1 < 0$, leading to vertical oscillations of the spot. That is, we observe that when $a_2 < 0$, the $(0,1)$ oscillation mode is dominant, consistent with the result of \S \ref{perturbeddisk}. We emphasize that the only difference between Figs. \ref{a2positive} and \ref{a2negative} is the sign of the $a_2$ coefficient in \eqref{pertf}. Furthermore, we note that $|a_2|$ is relatively small in comparison to the coefficients of the $\cos 3\theta$ and $\cos 4\theta$ terms, and yet it is the coefficient that dictates which mode of oscillation is preferred. This is due to the fact that the effect of the $\cos 2\theta$ perturbation enters at leading order in $\sigma$, while those of higher modes enter at $\mO(\eps^2)$.

\textbf{Experiment 2.} \, In left portion of Fig.~\ref{halfdisk1}, we show the unstable oscillations of one spot in the half disk when $\tau$ exceeds the Hopf threshold. On the right, we plot the $x_1$ and $x_2$ coordinates of the spot-center as a function of time;  observe that the initial oscillations at onset are only in the $x_1$-direction. Indeed, the predicted Hopf thresholds are $\hat{\tau} = 0.0712$ for the $(1,0)$ oscillatory mode and $\hat{\tau} = 0.1072$ for $(0,1)$ mode. The lower threshold for the $(1,0)$ mode indicates that it is the dominant mode of oscillation, which is what is observed numerically. The saturation of the oscillation amplitudes indicates that the Hopf bifurcation is supercritical, with the initial horizontal oscillations leading to a stable orbit with nonzero $x_1$ and $x_2$ components.

\begin{figure}[!ht]
     \centering
     \begin{subfigure}[b]{0.48\textwidth}
         \centering
         \includegraphics[width=\textwidth]{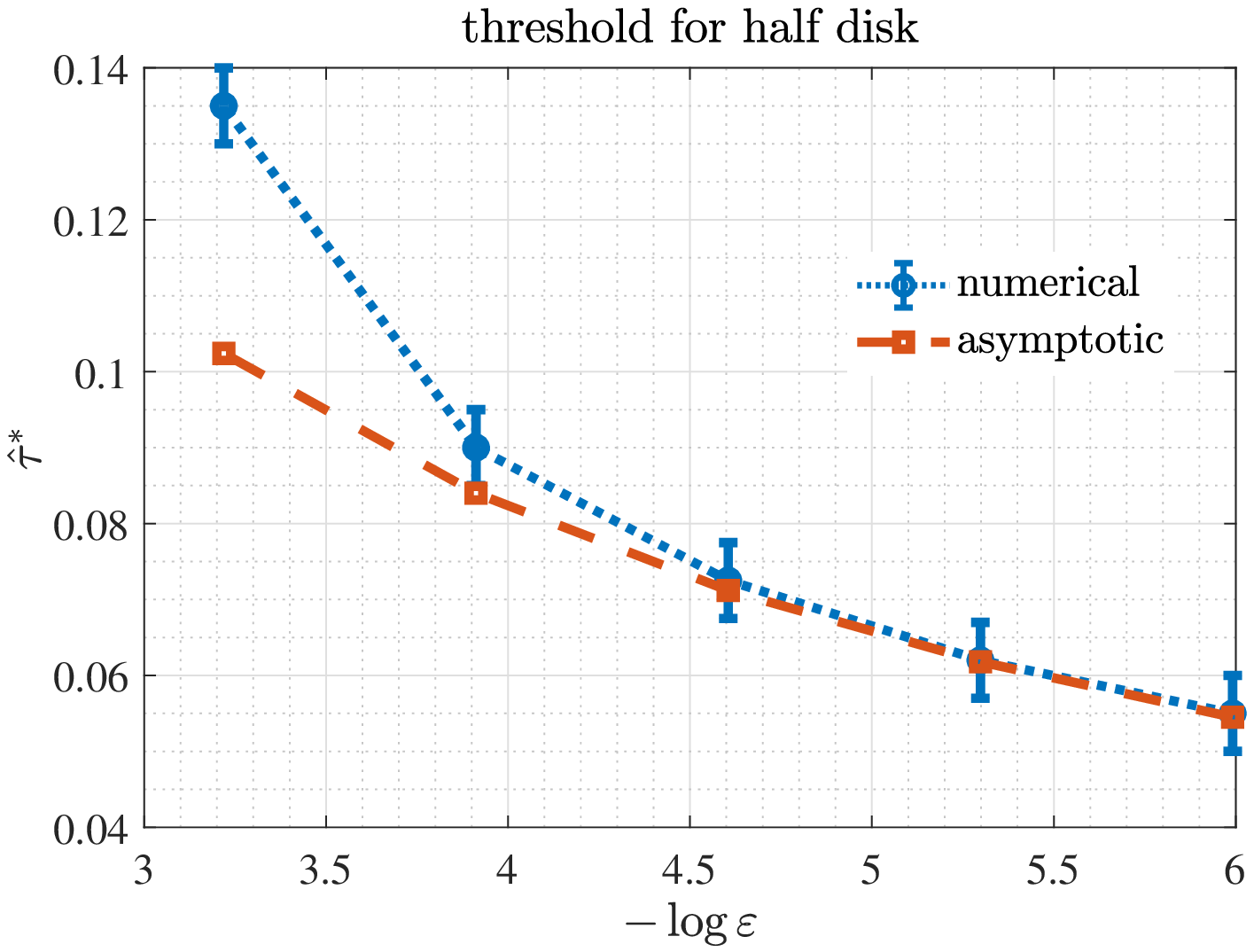}
         \caption{$\hat{\tau}^*$ vs. $-\log{\varepsilon}$}
         \label{tauvseps}
     \end{subfigure}
     \begin{subfigure}[b]{0.48\textwidth}
         \centering
         \includegraphics[width=\textwidth]{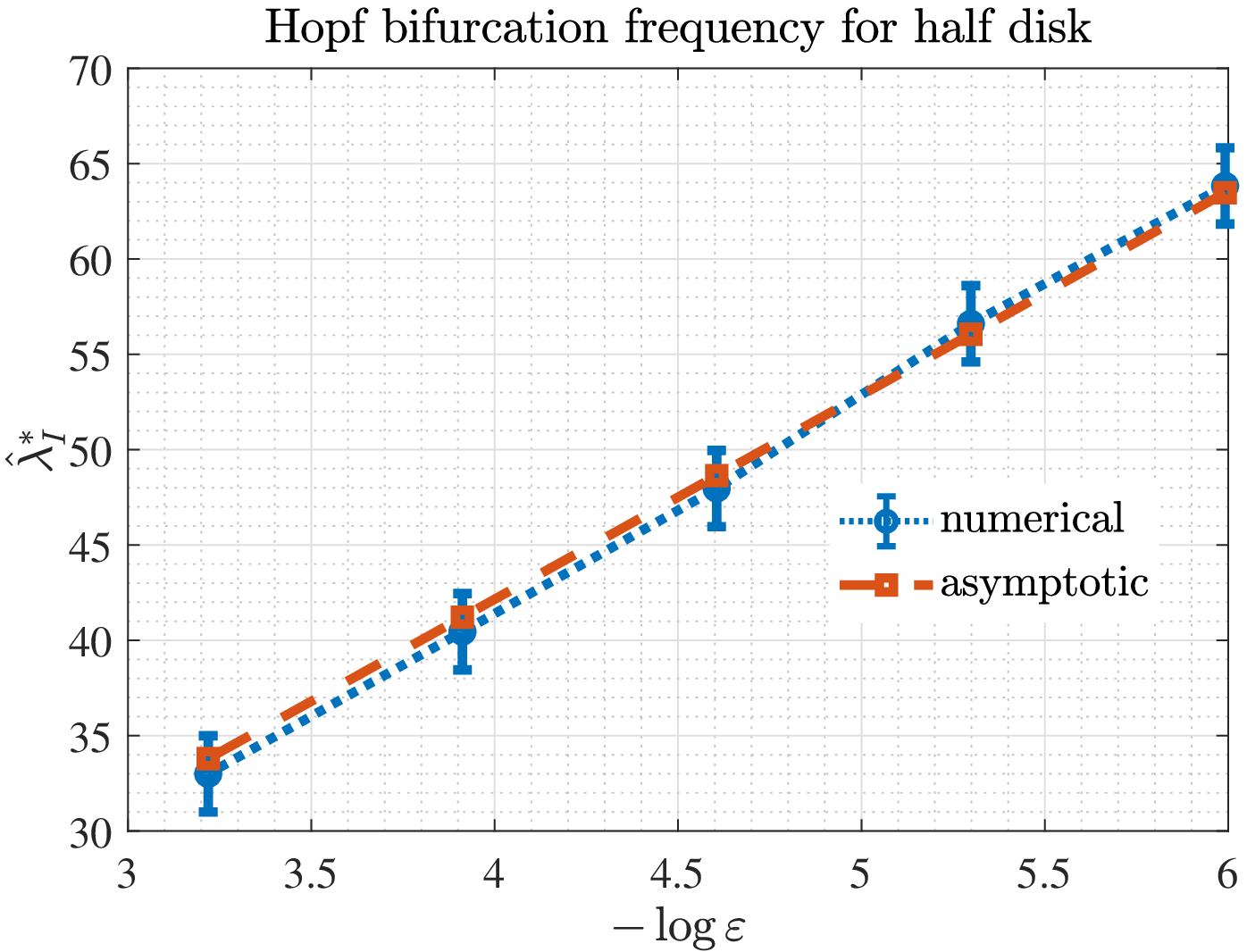}
         \caption{$\hat{\lambda}_I^*$ vs. $-\log{\varepsilon}$}
         \label{lambdavseps}
     \end{subfigure}
        \caption{Experiment 2 -- comparison of Hopf thresholds and corresponding frequency for $\varepsilon=0.04,~0.02,~0.01,~0.05,~0.025$. The blue error error bars indicate the range in which the exact threshold (obtained from PDE simulations) falls. The agreements in the Hopf stability threshold (a) as well as the corresponding frequency (b) between the predicted value and PDE simulations improve with decreasing $\varepsilon$.}
        \label{precision}
\end{figure}

\begin{figure}[!htb]
    \centering
    \includegraphics[width=0.98\textwidth]{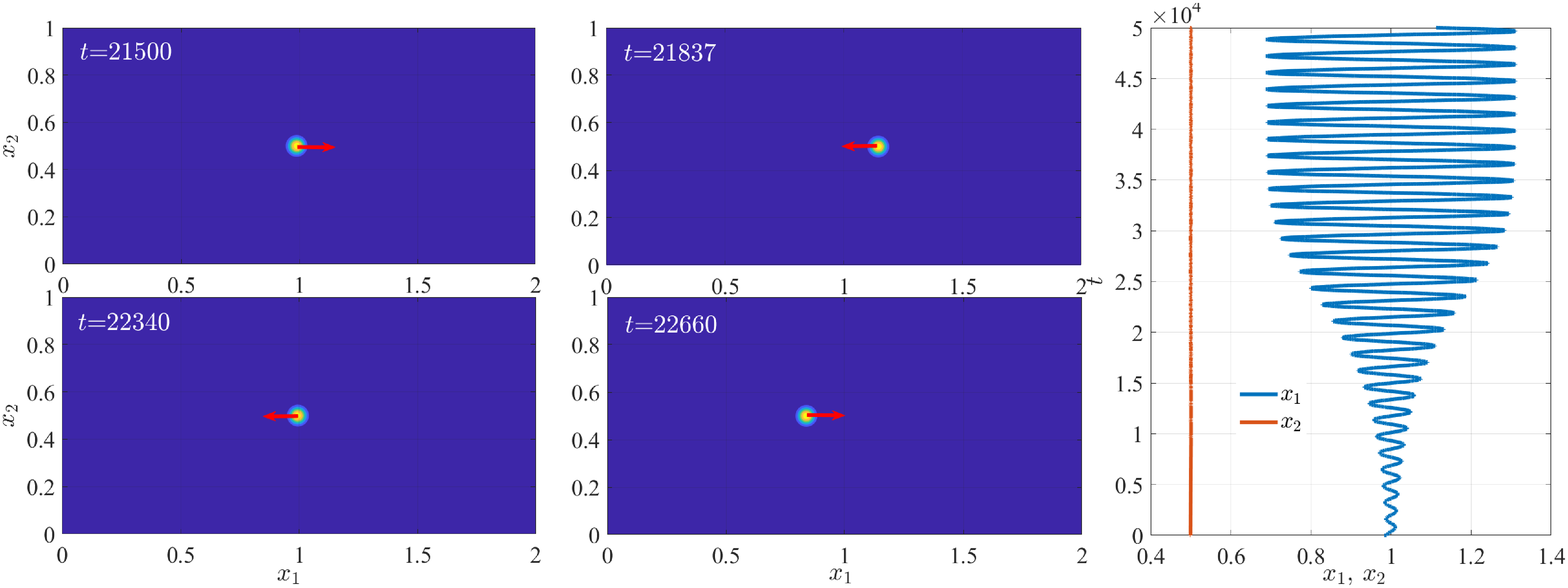}
        \caption{Experiment 3 -- numerical simulations performed for $\tau_0=0.063,~\varepsilon=0.01$ and $S= 4$. On the left, we show oscillatory motion of the spot center as $\tau_0$ exceeds the Hopf bifurcation threshold. Red arrows indicate the direction of motion. On the right, we plot the coordinates of the center of the spot, where the blue (red) curve is the $x_1$-coordinate ($x_2$-coordinate). The saturation of the oscillation amplitudes indicates that the Hopf bifurcation is supercritical, with the initial horizontal oscillations leading to a stable orbit with only a nonzero $x_1$ component.}
        \label{rectangle1}
\end{figure}

In Fig.~\ref{precision}, for the half-disk, we demonstrate the convergence with respect to $\varepsilon$ of the predicted values of the Hopf threshold $\hat{\tau}$ (Fig. \ref{tauvseps}) and corresponding bifurcation frequency $\lambda_I$ (Fig. \ref{lambdavseps}) to their exact values obtained from numerical computations. The blue error bars indicate the range in which the exact threshold falls. We note that all quantities plotted in Fig.~\ref{precision} have been scaled by $\eps^{-2}$; as such, the figures show that the error in the unscaled Hopf thresholds and frequencies scale as $\mO(\eps^2|\log\eps|)$ when $\eps$ is sufficiently small.

\textbf{Experiment 3.} \, In the left portion of Fig.~\ref{rectangle1}, we show Hopf oscillations of a spot in a rectangle of height $1$ and width $\ell = 2$. In our numerical simulations, we adjust the value of the feed-rate $A$ in \eqref{schnak} in inverse proportion to $\ell$ so as to keep the spot strength $S = (2\pi)^{-1} A|\Omega|$ constant. The dominant mode of oscillations is in the direction of the longer dimension, similar to what was shown for the perturbed disk of \S \ref{perturbeddisk}. On the right, we plot the $x_1$ and $x_2$ coordinates of the spot center. As was the case for the half-disk, the initial growth of amplitude oscillations saturates, suggesting that the Hopf bifurcation is supercritical. Due to the symmetry of the rectangle, the long-time oscillations occur only in the $x_1$ component. 

\begin{figure}[!htb]
    \centering
    \includegraphics[width=0.6\textwidth, height=0.35\textwidth]{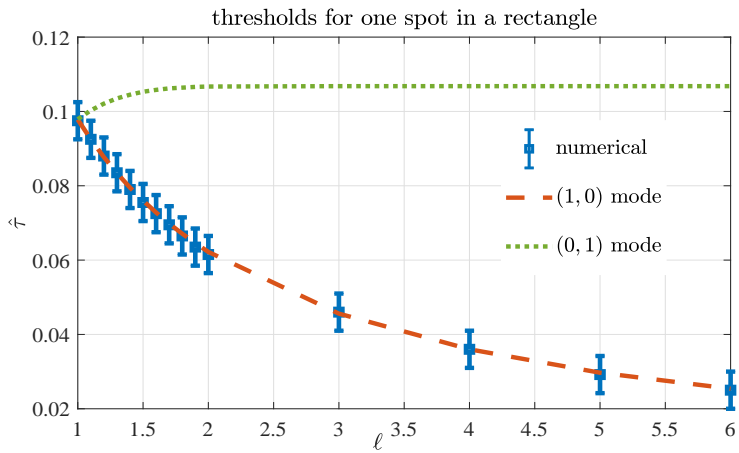}
    \caption{Experiment 3 -- thresholds for one spot in a rectangle of unit height and width $\ell$. When $\ell = 1$, the thresholds corresponding to the $(1,0)$ (horizontal) and $(0,1)$ (vertical) oscillation modes are equal due to symmetry. When $\ell$ increases, the dominant mode of oscillation is the one parallel to the longer edge of the rectangle.}
    \label{TauVsLen}
\end{figure}

\begin{figure}[!htb]
    \centering
    \includegraphics[width=0.98\textwidth]{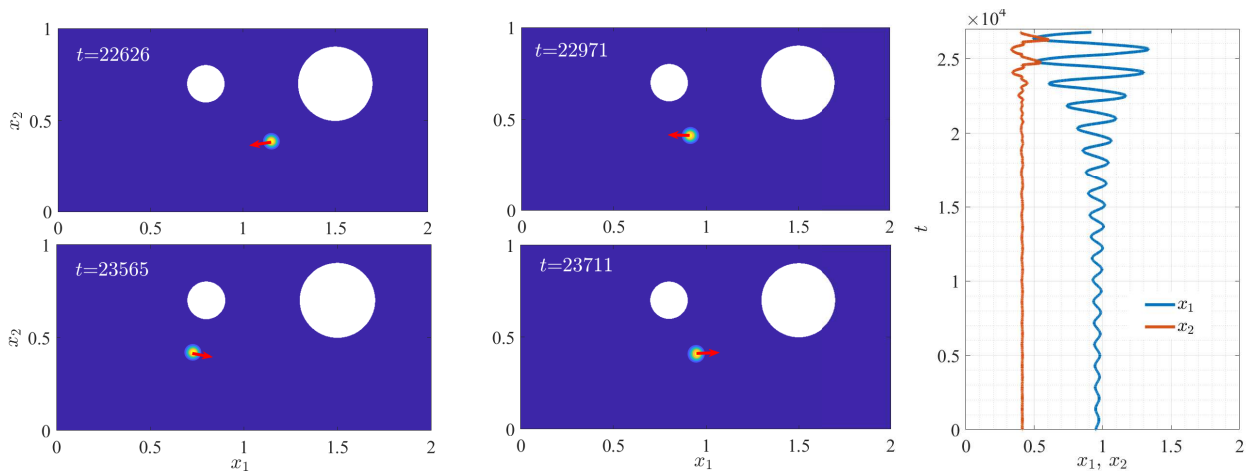}
        \caption{Experiment 4 -- numerical simulations performed for $\tau_0=0.063,~\varepsilon=0.01$ and $S= 4$. On the left, we show the oscillatory motion of the spot center as $\tau_0$ exceeds the Hopf bifurcation threshold. The red arrows indicate the direction of motion at the particular instant in time. On the right, we plot the coordinates of the spot center, where blue curve is the $x_1$-coordinate and red curve is the $x_2$-coordinate. The holes in the domain break the symmetry of the rectangle and cause nonzero long-time oscillations in both coordinates.}
        \label{rectangle2}
\end{figure}

In Fig.~\ref{TauVsLen}, we plot the predicted (red, green) and numerical (blue bars) Hopf bifurcation thresholds versus the length $\ell$ of the rectangle of unit height. The blue error error bars indicate the range in which the exact threshold falls as determined from PDE simulations. The two predicted thresholds correspond to the $(1,0)$ (horizontal) and $(0,1)$ (vertical) oscillation modes. When $\ell = 1$, the two thresholds are equal due to symmetry. As $\ell$ increases, the asymptotics results predict that the $(1,0)$ mode is the first to destabilize as $\tau_0$ is increased, in agreement with the oscillations parallel to the longer edge in the rectangle of Fig. \ref{rectangle1}.

We observe that both the Hopf threshold and the corresponding frequency decrease with increasing $\ell$. Indeed, as $\ell \to \infty$, we expect a zero eigenvalue corresponding to translational invariance in the $x_1$-direction. For this infinite strip, the only Hopf bifurcation is in the $x_2$-direction.

\textbf{Experiment 4.} \, In Fig. \ref{rectangle2}, we break the symmetry of the rectangle of Experiment 2 by removing two circular holes of different sizes. The feed-rate $A$ is set so that the strength of the spot is $S = 4$. The break in symmetry causes the spot center to shift away from the $(1, 0.5)$ location, and leads to long-time oscillations that are non-zero in both $x_1$- and $x_2$-components. The initial oscillations at onset, however, are still nearly horizontal, while the bifurcation threshold $\hat{\tau}^*$ is also very close to that of Experiment 3.

\subsection{Hopf bifurcation of multiple spots}
\begin{figure}[!htb]
     \centering
     \begin{subfigure}[b]{0.48\textwidth}
         \centering
         \includegraphics[width=\textwidth]{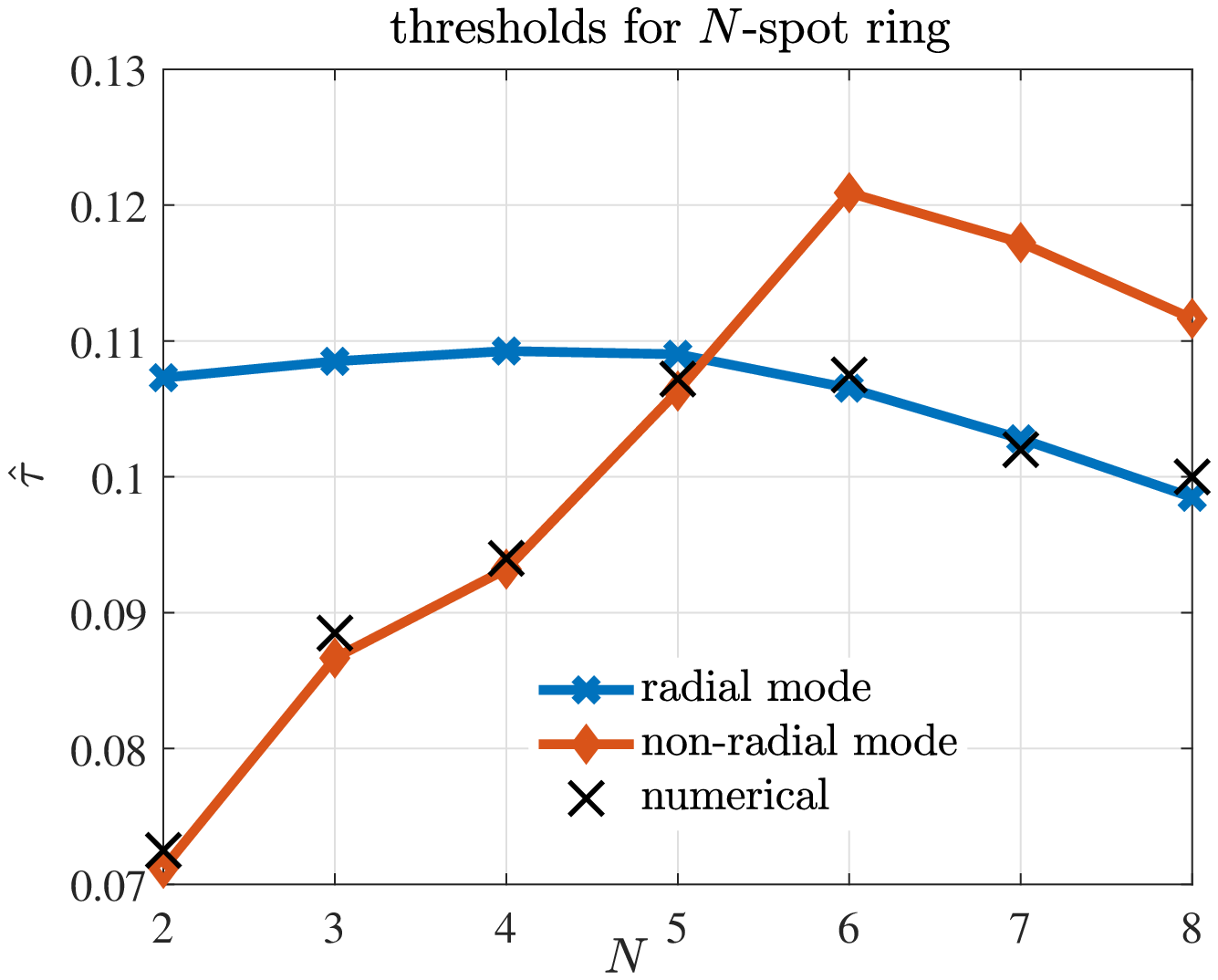}
         \caption{thresholds for $N$-spot ring in a disk}
         \label{crossingdisk}
     \end{subfigure}
     \begin{subfigure}[b]{0.48\textwidth}
         \centering
         \includegraphics[width=\textwidth]{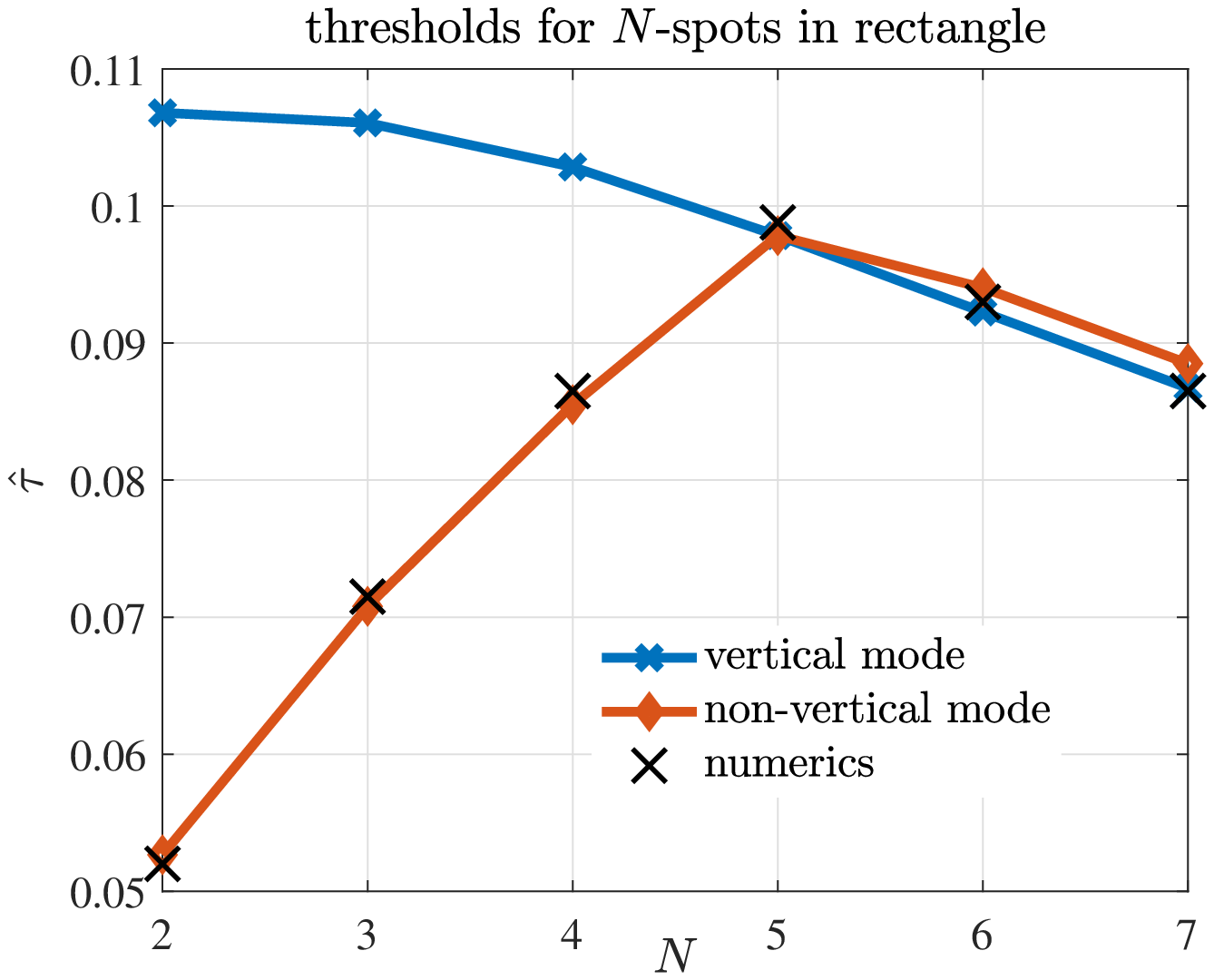}
         \caption{thresholds for $N$-spots in a rectangle}
         \label{crossingrect}
     \end{subfigure}
        \caption{Experiment 5 -- Hopf bifurcation thresholds for $N$-spots in a (a) unit disk and (b) rectangle. In (a), ``(near)-tangential'' refers to the mode in which spots oscillate in a direction tangent ($N$ even) or nearly tangent ($N$ odd) to the equilibrium ring, while ``next lowest found'' corresponds to the next lowest threshold (above that of the in-phase radial mode) that we were able to find in solving \eqref{finaleigblock}. The prediction that the in-phase radial mode is dominant for $N > 5$ is corroborated by the numerical simulations in Fig. \ref{UnitdiskM}. Similarly, in (b), the threshold plotted for the ``next lowest found'' is the lowest threshold we were able to find above that of the in-phase vertical mode. The emergence of the in-phase vertical mode as the dominant mode for $N \geq 6$ is observed in the numerical simulations of Fig. \ref{RectangleM}. In (a) and (b), the $\times$'s indicate the value of $\hat{\tau}_f$, the Hopf bifurcation value found from PDE simulations. }
        \label{crossing}
\end{figure}

\begin{figure}[!htb]
    \centering
    \includegraphics[width=0.95\textwidth]{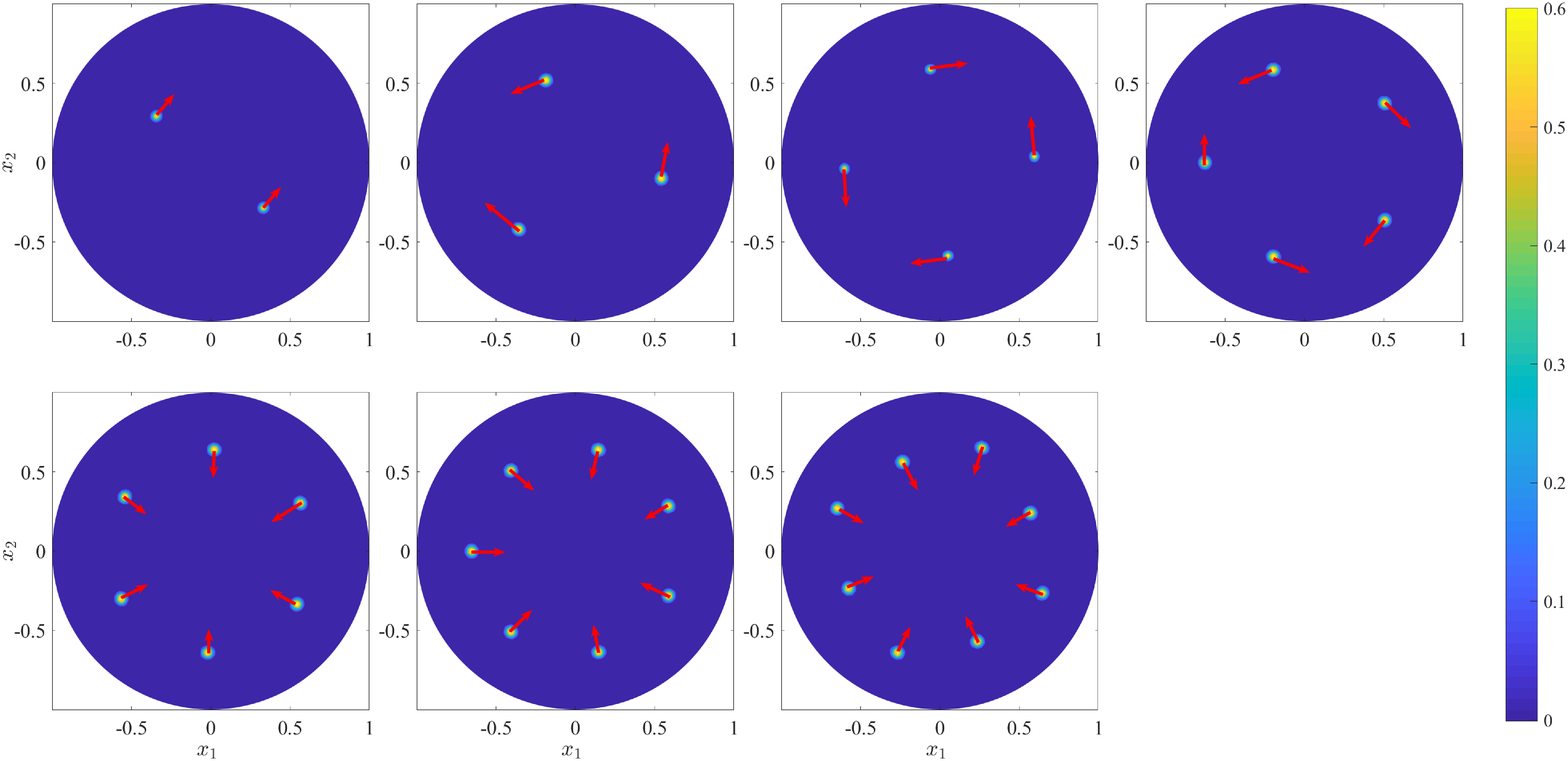}
    \caption{Experiment 5 -- spatial configurations of $N$-spots in a unit disk and their initial oscillation directions at onset for $\varepsilon=0.01,~S=4$. For $N=2,4$, the oscillations are along a direction tangential to the ring on which the spots are located. The $N = 3, 5$ configurations lack the symmetry to undergo perfectly tangential oscillations. The eigenvectors of \eqref{finaleigblock} are complex, indicating that the spots follow an elliptical orbit at onset. For $N=6,7,8$, the emergence of radial oscillations as the dominant mode for $N > 5$ is predicted by Fig. \ref{crossingdisk}. }
    \label{UnitdiskM}
\end{figure}

\begin{figure}[!htb]
    \centering
    \includegraphics[width=0.95\textwidth]{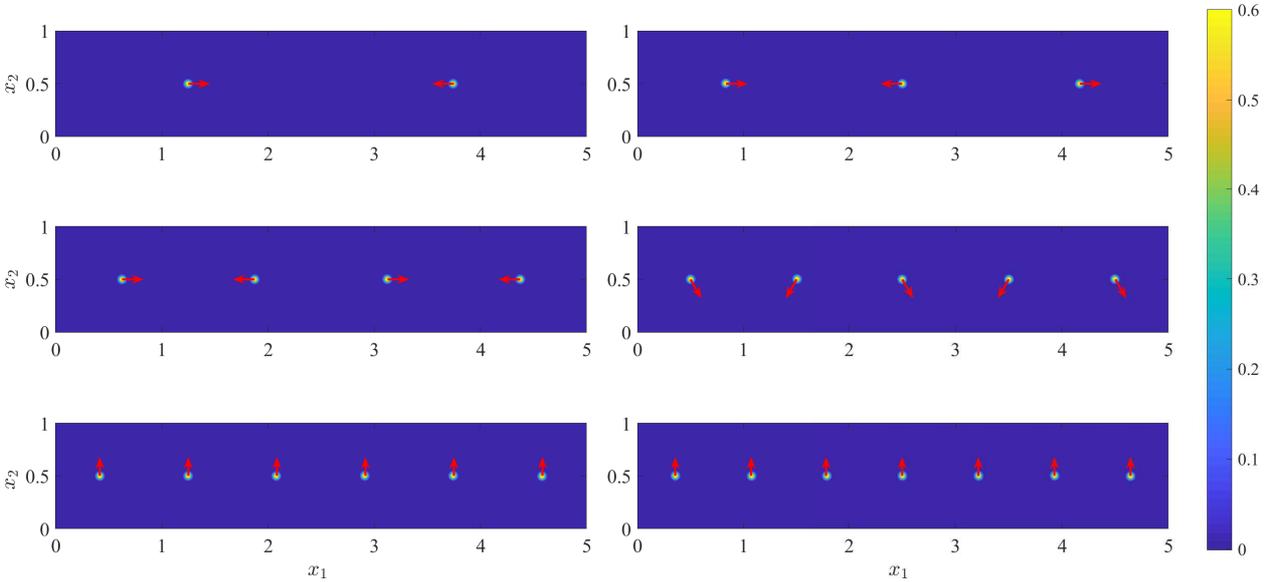}
    \caption{Experiment 5 -- spatial configurations and their initial oscillation directions at onset. There is a transition of the oscillation direction at $N=5$. For $N<5$, the oscillations are along the horizontal direction, while for $N>5$ the oscillations are following the vertical direction.  All profiles have one eigenfunction at the Hopf bifurcation except the case of 5-spot, where two eigenvectors are found, indicating that the initial oscillation have two directions. \red{The parameter values are $S=4$ and $\varepsilon=0.01$.}}
    \label{RectangleM}
\end{figure}

In this subsection, we investigate Hopf bifurcations of small eigenvalues of $N$-spot solutions to \eqref{schnak} along with the ensuing dynamics. While for a given $N$, multiple equilibrium configurations of spot locations and strengths are possible, we focus on the following configurations:

 \begin{enumerate}[resume]
    \item A ring of $N$ equally-spaced spots concentric with the unit disk (see Fig.~\ref{UnitdiskM}), and a line of $N$ equally-spaced spots along the width of a rectangle of height $1$ and width $5$ (see Fig. \ref{RectangleM}). Comparisons of stability thresholds are plotted in Fig. \ref{crossing} against $N$. We observe excellent agreement between asymptotic and numerical values.
    \item \label{RC-5H} Two spots in an asymmetric domain consisting of the same rectangle as that in Experiment 3 with five holes in the shape of disks. All boundary conditions are reflective. The holes create a barriers between the two spots. The predicted Hopf threshold is $\hat{\tau}^* = 0.0891$, while the threshold found from PDE simulations is $\hat{\tau}_f = 0.0887$.  
 \end{enumerate}

Spots in symmetric $N$-spot equilibrium configurations all have equal strength $S$ given by $S = \frac{A}{2N\pi|\Omega|}$. In the asymmetric domains of Experiment \ref{RC-5H}, spot strengths are generally unequal and must be determined by a numerical solution of the nonlinear system \eqref{Sj}, \eqref{eqbm1}, and \eqref{matchjj}. In the figures containing snapshots of solutions, the red arrows indicate the initial direction of oscillation at onset of instability. Arrows of different sizes indicate the relative amplitudes of oscillation between the different spots.

\textbf{Experiment 5.} \, In this experiment, we investigate dominant oscillation modes of symmetric $N$-spot equilibria. In particular, we show that our theory correctly predicts the switching of dominant modes as $N$ is increased. In Fig. \ref{crossingdisk}, we plot two bifurcation thresholds for an $N$-spot equilibrium arranged in a ring concentric with the unit disk (see Fig. \ref{UnitdiskM}). The blue curve corresponds to the radial mode of oscillations characterized by in-phase oscillations in the radial direction. For $N \leq 5$, the red diamonds denote the thresholds for the (near)-tangential mode in which spots oscillate (approximately) along the tangent of the equilibrium ring. As this threshold is lower than that of the radial mode, we expect that this mode emerges first as $\tau_0$ is increased. This is indeed observed in the first row of Fig. \ref{UnitdiskM}. The symmetry of even-numbered configurations allows this mode to be exactly tangential; for odd-numbered configurations, the spots undergo an elliptical orbit of high eccentricity. This elliptic orbit at onset is consistent with the fact the eigenvectors $\mathbf{a}^*$ of the eigenvalue problem \eqref{finaleigblock} for the $N = 3,5$ modes are complex with imaginary components small in comparison to the real components. Theoretical radii of these ring configurations were obtained from Eq. (16) of \cite{kolokolnikov2022ring}, and were used to initialize the simulations.

For $N > 5$, the purple squares of Fig. \ref{crossingdisk} denote the next lowest threshold (above that of the in-phase radial threshold) that we were able to find in the nonlinear system \eqref{det0}. For such $N$, the radial mode is the lower threshold, in which case we expect that this mode emerges first as $\tau_0$ is increased. This is seen in the second row of Fig. \ref{UnitdiskM}. The black $\times$'s in  Fig. \ref{crossingdisk} denote the value of $\tau_0$ at which the first Hopf bifurcation was encountered in our PDE simulations of \eqref{schnak}. The agreement, including the switch of dominance from the non-radial to radial mode of oscillation at $N = 6$, is excellent.

We observe a similar switch in mode-dominance for $N$-spots equally spaced along the center line parallel to the longer edge of the rectangle. In Fig. \ref{crossingrect}, the blue curve denotes the thresholds for the in-phase vertical oscillation mode in which all spots oscillate vertically in-phase (last row of Fig. \ref{RectangleM}). For $N \leq 5$, this threshold is preceded by a horizontal mode in which the spots oscillate horizontally and out-of-phase with their neighbors (red diamonds). For $N > 5$, we plot the next lowest threshold (above that of the in-phase vertical mode) that we were able to find in \eqref{finaleigblock} (purple squares). The black $\times$'s denote the value of $\tau_0$ at which the first Hopf bifurcation was encountered in our PDE simulations of \eqref{schnak}. We highlight the coincidence of the two thresholds at $N = 5$; indeed, when $N = 5$ in a rectangle of unit height and length $5$, the symmetry dictates that the two thresholds equal the $\ell = 1$ threshold of Fig. \ref{TauVsLen}. As such, the oscillations observed in the fourth image of Fig. \ref{RectangleM} is a linear combination of the $(1,0)$ and $(0,1)$ modes of a single spot in the unit square. 

This experiment illustrates what was concluded in \S \ref{perturbeddisk} for a single spot in a perturbed disk - the dominant mode of oscillation appears to be along the direction(s) in which there is more separation between spots or between a spot and the boundary. We also observe in these two scenarios that, when possible, the even mode was the dominant mode of oscillation. That is, each dominant mode can be replicated by a single spot in a correctly chosen domain with pure Neumann boundary conditions; i.e., a wedge of angle $2\pi/N$ for the $N$-spot ring in the unit disk, and a rectangle of unit height and length $\ell = 5/N$.

%

\begin{figure}[!htb]
     \centering
     \begin{subfigure}[b]{0.68\textwidth}
         \centering
         \includegraphics[width=\textwidth]{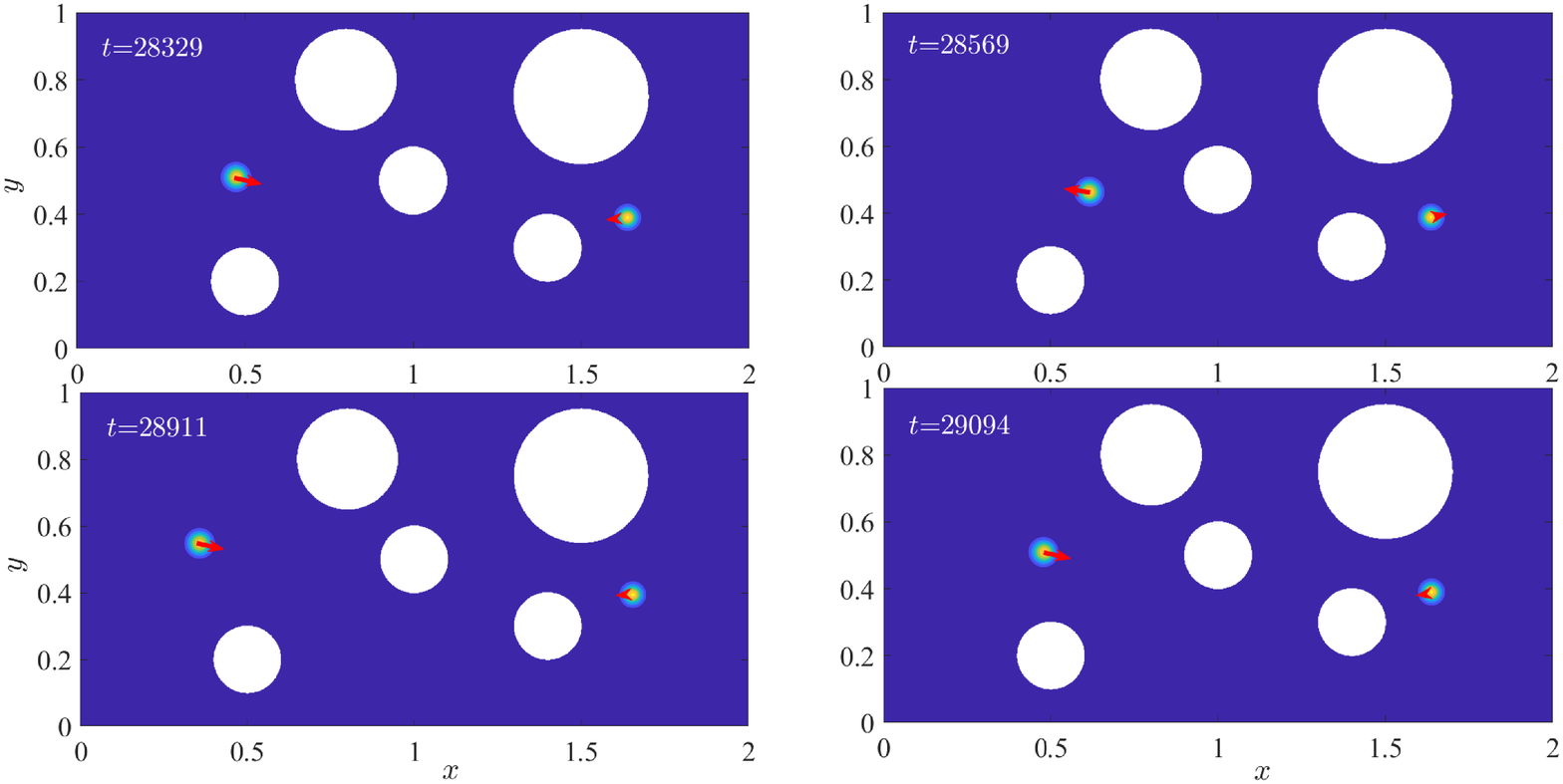}
         \label{rc5snapshots}
     \end{subfigure}
     \begin{subfigure}[b]{0.3\textwidth}
         \centering
         \includegraphics[width=\textwidth]{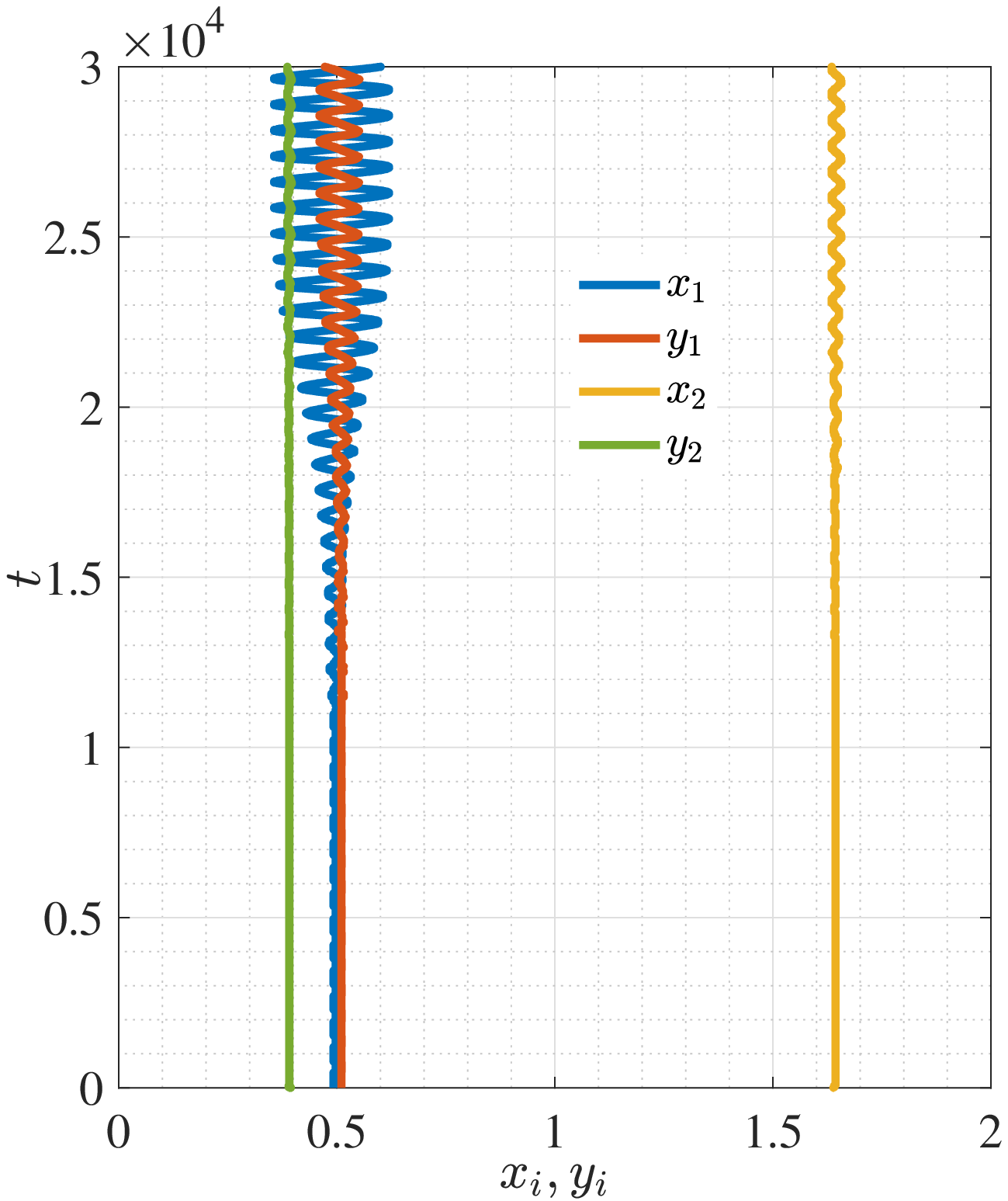}
         \label{rc5lefttrajectory}
     \end{subfigure}
        \caption{Experiment 6 -- numerical simulations performed for $\tau_0=0.09,~\varepsilon=0.01$. The strengths for the left spot and right spot are approximately $4$ and $3$, respectively. \red{On the left, we show snapshots of the two-spot profile at different times, along with the direction of motions indicated by the red arrows. The disparity in the sizes of the arrows convey the significantly smaller amplitude of oscillation of the right spot in comparison to that of the left. On the right, we plot the coordinates of the spot centers, where $(x_1,y_1)$ refers to the coordinates of the center of the left spot, while $(x_2,y_2)$ refers to those of the center of the right spot.} }
        \label{rectangleRc5holes}
\end{figure}

\textbf{Experiment 6.} \, In this experiment, we show the full generality of our stability result \eqref{finaleigblock} in which the Hessian terms of the Helmholtz and Green's functions were computed from a finite elements methods. Fig. \ref{rectangleRc5holes} shows two spots in a non-simply connected domain. While the threshold does not deviate significantly from that of one spot in a unit square, mode of oscillation as illustrated in the figure is rather different. The oscillation of the left spot is influenced by the orientation of the two nearest holes and has a vertical component as a result. The right spot is isolated in a smaller region and thus undergoes an oscillation of significantly smaller amplitude in comparison to the left spot (as indicated by the size of the arrows). A weak coupling still appears to be present, however, as the directions of oscillation still have remnants of the even mode of oscillation observed in the absence of the barriers.


\section{Discussion} \label{conc}

Through a formal asymptotic analysis, we have derived a $2N\times 2N$ complex matrix eigenvalue problem yielding the Hopf bifurcation threshold along with frequency and mode of oscillation for the slow oscillatory translation instabilities of $N$-spot equilibrium solutions to the Schnakenberg reaction-diffusion system with insulating boundary conditions. This result is valid for general, flat and bounded two-dimensional domains, and is a generalization of that derived for a single spot in a unit disk in \cite{xie2017moving}. This general result requires a more intricate analysis that accounts for domain asymmetries as well as spot-spot interactions.

For $N=1$, we showed that the matrix eigenvalue problem for bifurcation threshold and frequency reduces to the complex scalar problem derived in \cite{xie2017moving} for the unit disk. We then extended this analysis to that for a perturbed unit disk with radius in polar coordinates $r = 1 + \sigma f(\theta)$, $\theta \in [0,2\pi)$, $\eps \ll \sigma \ll 1$, where $f(\theta)$ is a $2\pi$-periodic function. We found that, at leading order in $\sigma$, only the coefficients $a_2 \equiv \pi^{-1}\int_0^{2\pi} \!f(\theta)\cos 2\theta \, d\theta$ and $b_2 \equiv \pi^{-1}\int_0^{2\pi} \!f(\theta)\sin 2\theta \, d\theta$  change the eigenvalue problem at leading order in $\sigma$. If $a_2 = b_2  = 0$, the bifurcation threshold will remain unchanged at leading order. When they are not both zero, the bifurcation threshold will decrease, and the preferred direction of oscillation will be along the line that makes an angle $\phi/2$ with the $x_1$-axis, where $\cos\phi = a_2/\sqrt{a_2^2+b_2^2}$ and $\sin\phi = b_2/\sqrt{a_2^2 + b_2^2}$. All of our asymptotic results were numerically confirmed from finite elements solutions of the full Schnakenberg PDE \eqref{schnak} using FlexPDE. 

Our analysis is valid for the Schnakenberg model with homogeneous feed-rate in a general flat $2$-dimensional domain of finite size. It would be interesting to investigate how various heterogeneous effects impact the stability threshold as well as the spot dynamics at onset. For example, \cite{wong2021spot} employs a hybrid asymptotic-numerical method to reveal novel dynamics and bifurcations of spot solutions when in the presence of a strongly localized feed-rate or small holes in the domain through which chemicals can leak. Heterogeneity can also come in the form of surface curvature. In this case, the integration of microlocal techniques into our asymptotic framework would be necessary to accurately compute Hessian terms of relevant Green's functions in the matrix eigenvalue problem. These techniques were first employed to compute linear terms of Green's function expansions for the purposes of predicting spot dynamics on curved surfaces in \cite{tzou2019spot} and \cite{tzou2020analysis}. Lastly, the small eigenvalues of spot clusters has not yet been analyzed. These clusters may form in the presence of a spatially dependent advection term in the PDE, or, as was analyzed in \cite{kolokolnikov2020hexagonal}, a spatially varying potential in the Gierer-Meinhardt model.

While we performed our analysis on the Schnakenberg model, a similar analysis would be possible for other activator-inhibitor reaction-diffusion models such as the Gray-Scott and Brusselator models. Furthermore, the hybrid method employed in this paper can be extended to compute small eigenvalues of $3$-dimensional spot patterns. For 3-D domains, the slow dynamics of quasi-equilibrium spot patterns along the stability and bifurcation structure of equilibrium configurations were analyzed for the Schnakenberg model in \cite{tzou2017stability} and more recently for the Gierer-Meinhardt model in \cite{gomez2021asymptotic}.

It would also be interesting to study the weakly nonlinear behavior of the spot dynamics beyond the linear stability regime. Such a theory has been developed in \cite{veerman2015breathing} for oscillatory amplitude instabilities in one-dimension. Numerical results of \cite{xie2017moving} showed the impact of the domain geometry on the long-time periodic orbits of a single spot when the Hopf bifurcation parameter $\tau$ was slightly above threshold. In particular, they observed primarily elliptical orbits for a single spot inside the unit disk, with a more exotic periodic orbit when the domain was a square. It would be an interesting analysis to characterize the periodic orbits that are admitted by a given domain and to investigate their stability. Another possible weakly nonlinear study could be performed near a codimension-2 point where the Hopf threshold for the small eigenvalues is equal to that for the large eigenvalues which lead to amplitude oscillations.

Lastly, it has been well-documented that equilibrium configurations of $N$-spot patterns in the Schnakenberg model correspond to (locally) optimal target configurations in the narrow escape optimization problem. See \cite{cheviakov2010asymptotic, ChevWard2010, tzou2017stability}. Furthermore, \cite{xie2017moving} and \cite{tzou2015mean} along with \cite{tzou2014first} establish connections between oscillatory translational instabilities of single spot patterns with a certain optimization problem of  a deterministically target on a one-dimensional interval and two-dimensional unit disk. It would be interesting to observe if this correspondence persists for multi-spot patterns and in more general two-dimensional domains.

\bibliographystyle{siam}
\bibliography{small_eig}

\end{document}